\documentclass[12pt,a4paper]{article}

\usepackage[left=2cm,right=2cm,    top=2cm,bottom=2cm,bindingoffset=0cm]{geometry}

\usepackage{subfig}
\usepackage{flushend}
\usepackage{indentfirst}
\usepackage{graphics}
\usepackage{amsmath}
\usepackage{epstopdf}
\usepackage{float}
\usepackage{amssymb}
\usepackage[font=footnotesize,labelfont=bf]{caption}
\usepackage[labelsep=period]{caption}

\begin{document}

\newtheorem{ozn}{Definition}[section]
\newtheorem{thm}{Theorem}[section]
\newtheorem{prop}{Proposition}[section]
\newtheorem{nas}{Corollary}[section]
\newtheorem{zau}{Remark}[section]
\newtheorem{lema}{Lemma}[section]
\newtheorem{pry}{Example}[section]

\newcommand{\eps}{\varepsilon}
\newcommand{\me}{\mathbf}
\newcommand{\mr}{\mathbb}
\newcommand{\mt}{\mathsf}
\newcommand{\md}{\mathcal}
\newcommand{\ld}{\left}
\newcommand{\rd}{\right}
\newcommand{\kd}{\bigtriangleup}
\newcommand{\ip}{\int_{-\pi}^{\pi}}
\newcommand{\iii}{\int_{-\infty}^{\infty}}

\newcommand{\be}{\begin{equation}}
\newcommand{\ee}{\end{equation}}

\newcommand{\bem}{\begin{multline}}
\newcommand{\eem}{\end{multline}}

\newcommand{\bml}{\begin{multline*}}
\newcommand{\eml}{\end{multline*}}

\newcommand{\beg}{\begin{gather}}
\newcommand{\eeg}{\end{gather}}

\title{\textbf{Robust interpolation of sequences with periodically stationary multiplicative seasonal increments}}

\date{}

\maketitle

\noindent  Carpathian Math. Publ. 2022, 14 (1), 105--126.\\
doi:10.15330/cmp.14.1.105-126

\vspace{20pt}

\author{\textbf{Maksym Luz}$^{1}$, \textbf{Mykhailo Moklyachuk}$^{2,*}$,\\\\
 {$^{1}$BNP Paribas Cardif, Kyiv, Ukraine, \\
$^{2,*}$Taras Shevchenko National University of Kyiv, Kyiv 01601, Ukraine}\\
$^{*}$Corresponding Author: Moklyachuk@gmail.com}\\\\\\

\noindent \textbf{Abstract.} \hspace{2pt}
We consider stochastic sequences with periodically stationary generalized multiple increments of fractional order which combines cyclostationary, multi-seasonal, integrated and fractionally integrated patterns. We solve the interpolation problem for linear functionals
constructed from unobserved values of a stochastic sequence of this type  based on observations of the sequence with a periodically stationary noise sequence.
For sequences with known matrices of spectral densities, we obtain formulas for calculating values of the mean square errors and the spectral characteristics of the optimal interpolation of the functionals.
Formulas that determine the least favorable spectral densities and the minimax (robust) spectral
characteristics of the optimal linear interpolation of the functionals
are proposed in the case where spectral densities of the sequences
are not exactly known while some sets of admissible spectral densities are given.\\

\noindent \textbf{Keywords:} \hspace{2pt} periodically stationary increments, SARFIMA, fractional integration, filtering, optimal linear
estimate, mean square error, least favourable spectral density
matrix, minimax spectral characteristic\\

\noindent \textbf{AMS Subject Classification:} \hspace{2pt}
 Primary: 60G10, 60G25, 60G35, Secondary: 62M20, 93E10

\maketitle

\vspace{10pt}

\noindent {\bf{DOI}}:10.15330/cmp.14.1.105-126


\section*{Introduction}

Non-stationary and long memory time series models are wildly used in different fields of economics, finance, climatology, air pollution, signal processing etc.
(see, for example, papers by Dudek and Hurd \cite{Dudek-Hurd}, Johansen and Nielsen  \cite{Johansen}, Reisen et al.\cite{Reisen2014}).
A core example -- a general multiplicative model, or $SARIMA (p, d, q)\times(P, D, Q)_s$ -- was introduced in the book by Box and Jenkins \cite{Box_Jenkins}. It includes both  integrated and  seasonal factors:
\be
 \Psi (B ^ s) \psi (B) (1-B) ^ d (1-B ^ s) ^ Dx_t = \Theta (B ^ s) \theta (B) \eps_t, \label {seasonal_3_model}
 \ee
 where $ \eps_t $, $ t \in \mr Z $ is a sequence of zero mean i.i.d. random variables,
and where $ \Psi (z) $ and $ \Theta (z) $ are two polynomials of degrees of $ P $ and $ Q $ respectively which have roots outside the unit circle.
The parameters $d$ and $D$ are allowed to be fractional. In the case where $|d+D|<1/2$ and $|D|<1/2$, the process  $\eqref{seasonal_3_model}$ is stationary and invertible.
The paper  by Porter-Hudak \cite{Porter-Hudak} illustrates an application of a seasonal ARFIMA model to the analysis of the monetary aggregates used by U.S. Federal Reserve.
Another model of fractional integration is GARMA processes  described by the equation (see Gray,  Cheng and  Woodward \cite{Gray})
\be
 (1-2uB+B^2) ^ dx_t = \eps_t,\quad |u|\leq1. \label{GARMA_model} \ee

For the resent results dedicated to the statistical inference for seasonal long-memory sequences,
we refer  to the paper by Tsai,  Rachinger and  Lin \cite{Tsai}, who developed methods of estimation of parameters in case of measurement errors.
In their paper Baillie,   Kongcharoen and  Kapetanios \cite{Baillie} compared MLE and semiparametric estimation procedures for prediction problems based on ARFIMA models.
 Based on simulation study, they indicate better performance of  MLE predictor than the one  based on the two-step local Whittle estimation.
Hassler and  Pohle \cite{Hassler} (see also Hassler \cite{Hassler_book}) assess a predictive  performance of various methods of  forecasting  of inflation and return volatility time series and show strong evidences for models with a fractional integration component.

Another type of non-stationary processes are periodically correlated, or cyclostationary, processes introduced by Gladyshev \cite{Gladyshev}, which belong to the class of processes with time-dependent spectrum and are widely used in signal processing and communications (see Napolitano \cite{Napolitano} for a review of recent works on cyclostationarity and its applications).
Periodic time series are also considered as extension of seasonal models \cite{Baek,Basawa,Lund,Osborn}.

One of the fields of interests related to time series analysis is optimal filtering. It aims to remove the unobserved components, such as trends, seasonality or noise signal, from the observed data  \cite{Ansley,Eiurridge}.

Methods of parameters estimations and filtering usually do not take into account the issues arising from real data, namely, the presence of outliers, measurement errors, incomplete information about the spectral, or model, structure etc. From this point of view, we see an increasing interest to robust methods of estimation that are reasonable in such cases (see
Reisen,  et al. \cite{Reisen2018}, Solci at al. \cite{Solci} for the examples of robust estimates of SARIMA and PAR models). Grenander \cite{Grenander},
 Hosoya \cite{Hosoya}, Kassam  \cite{Kas1982}, Franke \cite{Franke1985}, Vastola and  Poor \cite{VastPoor1984}, Moklyachuk \cite{Moklyachuk,Moklyachuk2015},
Luz and Moklyachuk \cite{Luz_Mokl_filt3},
 Liu et al. \cite{Liu} studied minimax  extrapolation, interpolation and filtering problems for stationary sequences and processes.

This article is dedicated to the robust interpolation of stochastic sequences with periodically stationary long memory multiple seasonal increments, or sequences with periodically stationary general multiplicative (GM) increments, introduced by Luz and Moklyachuk \cite{Luz_Mokl_extra_GMI}. In the recent years, models with multiple seasonal and periodic patterns are in scope of interest, see  Dudek \cite{Dudek}, Gould et al. \cite{Gould} and Reisen et al. \cite{Reisen2014}, Hurd and Piparas \cite{Hurd}. This  motivates us to study the time series combining a periodic structure of the covariation function as well as multiple seasonal factors. The discussed problem is a natural continuation of the researches on minimax interpolation and filtering of stationary vector-valued processes, periodically correlated processes and processes with stationary increments have been performed by  Moklyachuk and Masyutka \cite{Mokl_Mas_filt},  Moklyachuk and  Golichenko (Dubovetska) \cite{Dubovetska_filt},  Luz and Moklyachuk \cite{Luz_Mokl_filt,Luz_Mokl_book} respectively. We also mention the works by Moklyachuk, Masyutka and Sidei \cite{Sidei_filt,Sidei_book,Nas_Mokl_Sidei},
 who derive minimax estimates of stationary processes from observations with missed values, and the work  by Moklyachuk and  Kozak \cite{Kozak_Mokl}, who have studied interpolation problem for stochastic sequences with periodically stationary increments.

The article is organized as follows.
In Section $\ref{spectral_ theory}$, we recall definitions of
 generalized multiple (GM)  increment sequence $\chi_{\overline{\mu},\overline{s}}^{(d)}(\vec{\xi}(m))$ and stochastic
sequences $\xi(m)$ with periodically stationary (periodically correlated, cyclostationary) GM increments.
The spectral theory of vector-valued GM increment sequences is discussed.
Section $\ref{classical_filtering}$ deals with the classical interpolation  problem for the linear functional  $A_N\xi $ which is constructed from unobserved values of the sequence $\xi(m)$ when the spectral densities of the sequence $\xi(m)$ and a noise  sequence $\eta(m)$ are known.
Estimates are obtained by applying the Hilbert space projection technique to the vector-valued sequence $\vec \xi(m)+ \vec \eta(m)$ with stationary GM  increments under the stationary noise sequence $\vec \eta(m)$ uncorrelated with $\vec \xi(m)$.
The case of non-stationary fractional integration is discussed as well.
Section $\ref{minimax_filtering}$ is dedicated to the minimax (robust) estimates in cases, where spectral densities of sequences are not exactly known
while some sets of admissible spectral densities are specified. We illustrate the proposed technique on the particular types of the sets, which are generalizations of the sets of admissible spectral densities described in a survey article by
Kassam and Poor \cite{Kassam_Poor} for stationary stochastic processes..

\section{Stochastic sequences with periodically stationary generalized multiple increments}\label{spectral_ theory}

\subsection{Preliminary notations and definitions}
 Consider a   stochastic sequence $\xi(m)$, $m\in\mathbb Z$ defined on a probability space $(\Omega, \cal F, \mathbb P)$.
 Denote by $B_{\mu}$ a backward shift operator   with the step $\mu\in
\mathbb Z$, such that $B_{\mu}\xi(m)=\xi(m-\mu)$; $B:=B_1$. Then $B_{\mu}^s=B_{\mu}B_{\mu}\cdot\ldots\cdot B_{\mu}$.
Define an incremental operator $$\chi_{\overline{\mu},\overline{s}}^{(d)}(B)
=(1-B_{\mu_1}^{s_1})^{d_1}(1-B_{\mu_2}^{s_2})^{d_2}\cdots(1-B_{\mu_r}^{s_r})^{d_r},$$ where
$d:=d_1+d_2+\ldots+d_r$, $\overline{d}=(d_1,d_2,\ldots,d_r)\in (\mr N^*)^r$,
 $\overline{s}=(s_1,s_2,\ldots,s_r)\in (\mr N^*)^r$
and $\overline{\mu}=(\mu_1,\mu_2,\ldots,\mu_r)\in (\mr N^*)^r$ or $\in (\mr Z\setminus\mr N)^r$. Here $\mr N^*=\mr N\setminus\{0\}$.

Within the article, $\delta_{lp}$ denotes the Kronecker symbols, and ${n \choose l}=\frac{n!}{l!(n-l)!}$.

\begin{ozn}\label{def_multiplicative_Pryrist}
For a stochastic sequence $\xi(m)$, $m\in\mathbb Z$, the
sequence
\begin{eqnarray}
\nonumber
\chi_{\overline{\mu},\overline{s}}^{(d)}(\xi(m))&:=&\chi_{\overline{\mu},\overline{s}}^{(d)}(B)\xi(m)
=(1-B_{\mu_1}^{s_1})^{d_1}(1-B_{\mu_2}^{s_2})^{d_2}\cdots(1-B_{\mu_r}^{s_r})^{d_r}\xi(m)
\\&=&\sum_{l_1=0}^{d_1}\ldots \sum_{l_r=0}^{d_r}(-1)^{l_1+\ldots+ l_r}{d_1 \choose l_1}\cdots{d_r \choose l_r}\xi(m-\mu_1s_1l_1-\cdots-\mu_rs_rl_r)
\label{GM_Pryrist}
\end{eqnarray}
is called a \emph{stochastic  generalized multiple (GM)  increment sequence} of differentiation   order
$d$
with a fixed seasonal  vector $\overline{s}\in (\mr N^*)^r$
and a varying step $\overline{\mu}\in (\mr N^*)^r$ or $\in (\mr Z\setminus\mr N)^r$.
\end{ozn}

 The  multiplicative increment operator $\chi_{\overline{\mu},\overline{s}}^{(d)}(B)$ admits a representation
\[
\chi_{\overline{\mu},\overline{s}}^{(d)}(B)
=\prod_{i=1}^r(1-B_{\mu_i}^{s_i})^{d_i}
=\sum_{k=0}^{n(\gamma)}e_{\gamma}(k)B^k,
\]
where $n(\gamma):=\sum_{i=1}^r\mu_is_id_i$. The explicit representation of the coefficients $e_{\gamma}(k)$ is given in \cite{Luz_Mokl_extra_GMI}.

\begin{ozn}
\label{oznStPryrostu}
A stochastic GM increment sequence $\chi_{\overline{\mu},\overline{s}}^{(d)}(\xi(m))$  is called   a wide sense
stationary if the mathematical expectations
\begin{eqnarray*}
\mt E\chi_{\overline{\mu},\overline{s}}^{(d)}(\xi(m_0))& = &c^{(d)}_{\overline{s}}(\overline{\mu}),
\\
\mt E\chi_{\overline{\mu}_1,\overline{s}}^{(d)}(\xi(m_0+m))\chi_{\overline{\mu}_2,\overline{s}}^{(d)}(\xi(m_0))
& = & D^{(d)}_{\overline{s}}(m;\overline{\mu}_1,\overline{\mu}_2)
\end{eqnarray*}
exist for all $m_0,m,\overline{\mu},\overline{\mu}_1,\overline{\mu}_2$ and do not depend on $m_0$.
The function $c^{(d)}_{\overline{s}}(\overline{\mu})$ is called a mean value  and the function $D^{(d)}_{\overline{s}}(m;\overline{\mu}_1,\overline{\mu}_2)$ is
called a structural function of the stationary GM increment sequence (of a stochastic sequence with stationary GM increments).
\\
The stochastic sequence $\xi(m)$, $m\in\mathbb   Z$
determining the stationary GM increment sequence
$\chi_{\overline{\mu},\overline{s}}^{(d)}(\xi(m))$ by   \eqref{GM_Pryrist} is called a stochastic
sequence with stationary GM increments (or GM increment sequence of order $d$).
\end{ozn}

\begin{zau}
For spectral properties of one-pattern increment sequence $\chi_{\mu,1}^{(n)}(\xi(m)):=\xi^{(n)}(m,\mu)=(1-B_{\mu})^n\xi(m)$
see, e.g., \cite{Luz_Mokl_book}, p. 1-8; \cite{Gihman_Skorohod}, p. 48--60, 261--268; \cite{Yaglom}, p. 390--430.
The corresponding results for continuous time increment process $\xi^{(n)}(t,\tau)=(1-B_{\tau})^n\xi(t)$ are described in \cite{Yaglom:1955}, \cite{Yaglom}.
\end{zau}

\subsection{Definition and spectral representation of stochastic sequence with periodically stationary GM increment}

In this section, we present definition, justification and a brief review of the spectral theory of stochastic sequences with periodically stationary multiple seasonal increments.

\begin{ozn}
\label{OznPeriodProc}
A stochastic sequence $\xi(m)$, $m\in\mathbb Z$ is called a stochastic sequence with periodically stationary (periodically correlated) GM increments with period $T$ if the mathematical expectations
\begin{eqnarray*}
\mt E\chi_{\overline{\mu},T\overline{s}}^{(d)}(\xi(m+T)) & = & \mt E\chi_{\overline{\mu},T\overline{s}}^{(d)}(\xi(m))=c^{(d)}_{T\overline{s}}(m,\overline{\mu}),
\\
\mt E\chi_{\overline{\mu}_1,T\overline{s}}^{(d)}(\xi(m+T))\chi_{\overline{\mu}_2,T\overline{s}}^{(d)}(\xi(k+T))
& = & D^{(d)}_{T\overline{s}}(m+T,k+T;\overline{\mu}_1,\overline{\mu}_2)
= D^{(d)}_{T\overline{s}}(m,k;\overline{\mu}_1,\overline{\mu}_2)
\end{eqnarray*}
exist for every  $m,k,\overline{\mu}_1,\overline{\mu}_2$ and  $T>0$ is the least integer for which these equalities hold.
\end{ozn}

It follows from  Definition \ref{OznPeriodProc} that the sequence
\begin{equation}
\label{PerehidXi}
\xi_{p}(m)=\xi(mT+p-1), \quad p=1,2,\dots,T; \quad m\in\mathbb Z
\end{equation}
forms a vector-valued sequence
$\vec{\xi}(m)=\left\{\xi_{p}(m)\right\}_{p=1,2,\dots,T}, m\in\mathbb Z$
with stationary GM increments as follows:
\begin{eqnarray*}
\chi_{\overline{\mu},\overline{s}}^{(d)}(\xi_p(m))&=&\sum_{l_1=0}^{d_1}\ldots \sum_{l_r=0}^{d_r}(-1)^{l_1+\ldots+ l_r}{d_1 \choose l_1}\cdots{d_r \choose l_r}\xi_p(m-\mu_1s_1l_1-\ldots-\mu_rs_rl_r)
\\
&=&\sum_{l_1=0}^{d_1}\ldots \sum_{l_r=0}^{d_r}(-1)^{l_1+\ldots+ l_r}{d_1 \choose l_1}\cdots{d_r \choose l_r}\xi((m-\mu_1s_1l_1-\ldots-\mu_rs_rl_r)T+p-1)
\\
&=&\chi_{\overline{\mu},T\overline{s}}^{(d)}(\xi(mT+p-1)),\quad p=1,2,\dots,T,
\end{eqnarray*}
where $\chi_{\overline{\mu},\overline{s}}^{(d)}(\xi_p(m))$ is the GM increment of the $p$-th component of the vector-valued sequence $\vec{\xi}(m)$.

The following theorem describes the spectral structure of the GM increment  \cite{Karhunen}, \cite{Luz_Mokl_extra_GMI}.

\begin{thm}\label{thm1}
1. The mean value and structural function
 of the vector-valued stochastic stationary
GM increment sequence $\chi_{\overline{\mu},\overline{s}}^{(d)}(\vec{\xi}(m))$ can be represented in the form
\begin{eqnarray}
\label{serFnaR_vec}
c^{(d)}_{ \overline{s}}(\overline{\mu})& = &c\prod_{i=1}^r\mu_i^{d_i},
\\
\label{strFnaR_vec}
 D^{(d)}_{\overline{s}}(m;\overline{\mu}_1,\overline{\mu}_2)& = &\int_{-\pi}^{\pi}e^{i\lambda
m} \chi_{\overline{\mu}_1}^{(d)}(e^{-i\lambda})\chi_{\overline{\mu}_2}^{(d)}(e^{i\lambda})\frac{1}
{|\beta^{(d)}(i\lambda)|^2}dF(\lambda),
\end{eqnarray}
where
\[\chi_{\overline{\mu}}^{(d)}(e^{-i\lambda})=\prod_{j=1}^r(1-e^{-i\lambda\mu_js_j})^{d_j}, \quad \beta^{(d)}(i\lambda)= \prod_{j=1}^r\prod_{k_j=-[s_j/2]}^{[s_j/2]}(i\lambda-2\pi i k_j/s_j)^{d_j},
\]
 $c$ is a vector, $F(\lambda)$ is the matrix-valued spectral function of the stationary stochastic sequence $\chi_{\overline{\mu},\overline{s}}^{(d)}(\vec{\xi}(m))$. The vector $c$
and the matrix-valued function $F(\lambda)$ are determined uniquely by the GM
increment sequence $ \chi_{\overline{\mu},\overline{s}}^{(d)}(\vec \xi(m))$.

2. The stationary GM increment sequence $\chi_{\overline{\mu},\overline{s}}^{(d)}(\vec{\xi}(m))$ admits the spectral representation
\begin{equation}
\label{SpectrPred_vec}
\chi_{\overline{\mu},\overline{s}}^{(d)}(\vec{\xi}(m))
=\int_{-\pi}^{\pi}e^{im\lambda}\chi_{\overline{\mu}}^{(d)}(e^{-i\lambda})\frac{1}{\beta^{(d)}(i\lambda)}d\vec{Z}_{\xi^{(d)}}(\lambda),
\end{equation}
where $d\vec{Z}_{\xi^{(d)}}(\lambda)=\{Z_{ p}(\lambda)\}_{p=1}^{T}$ is a (vector-valued) stochastic process with uncorrelated increments on $[-\pi,\pi)$ connected with the spectral function $F(\lambda)$ by
the relation
\[
 \mt E(Z_{p}(\lambda_2)-Z_{p}(\lambda_1))(\overline{ Z_{q}(\lambda_2)-Z_{q}(\lambda_1)})
 =F_{pq}(\lambda_2)-F_{pq}(\lambda_1),\]
 \[  -\pi\leq \lambda_1<\lambda_2<\pi,\quad p,q=1,2,\dots,T.
\]
\end{thm}

Consider another vector-valued stochastic sequence with the stationary GM
increments $\vec \zeta (m)=\vec \xi(m)+\vec \eta(m)$, where $\vec\eta(m)$ is a vector-valued stationary stochastic sequence, uncorrelated with $\vec\xi(m)$, with a spectral representation
\[
 \vec\eta(m)=\int_{-\pi}^{\pi}e^{i\lambda m}dZ_{\eta}(\lambda),\]
 $Z_{\eta}(\lambda)=\{Z_{\eta,p}(\lambda)\}_{p=1}^T$, $\lambda\in [-\pi,\pi)$, is a stochastic process with uncorrelated increments, that corresponds to the spectral function $G(\lambda)$ \cite{Hannan}.
The stochastic stationary GM increment $\chi_{\overline{\mu},\overline{s}}^{(d)}(\vec{\zeta}(m))$ allows the spectral representation
\begin{eqnarray*}
 \chi_{\overline{\mu},\overline{s}}^{(d)}(\vec{\zeta}(m))&=&\int_{-\pi}^{\pi}e^{i\lambda m}\frac{\chi_{\overline{\mu}}^{(d)}(e^{-i\lambda})}{\beta^{(d)}(i\lambda)}
 dZ_{\xi^{(n)}}(\lambda)
  +\int_{-\pi}^{\pi}e^{i\lambda m}\chi_{\overline{\mu}}^{(d)}(e^{-i\lambda}) dZ_{\eta }(\lambda),
 \end{eqnarray*}
while $dZ_{\eta }(\lambda)=(\beta^{(d)}(i\lambda))^{-1} dZ_{\eta^{(n)}}(\lambda)$,
$\lambda\in[-\pi,\pi)$. Therefore, in the case where the spectral functions $F(\lambda)$ and $G(\lambda)$ have the spectral densities $f(\lambda)$ and $g(\lambda)$, the spectral density $p(\lambda)=\{p_{ij}(\lambda)\}_{i,j=1}^{T}$ of the stochastic sequence $\vec \zeta(m)$ is determined by the formula
\[
 p(\lambda)=f(\lambda)+|\beta^{(d)}(i\lambda)|^2g(\lambda).\]

\subsection{Sequences with GM fractional increments}\label{fractional_extrapolation}

In the previous subsection, we describe    the  GM increment sequence $\chi_{\overline{\mu},\overline{s}}^{(d)}(\vec{\xi}(m))$ of  the positive integer orders $(d_1,\ldots,d_r)$.
Here we consider the case of possibly fractional increment orders $d_i$.

Within the subsection, we put the step $\overline{\mu}=(1,1,\ldots,1)$. Represent the increment operator $\chi_{\overline{s}}^{(d)}(B)$  in the form
\be\label{FM_increment}
\chi_{\overline{s}}^{(R+D)}(B)=(1-B)^{R_0+D_0}\prod_{j=1}^r(1-B^{s_j})^{R_j+D_j},
\ee
where $(1-B)^{R_0+D_0}$ is an integrating component, $R_j$, $j=0,1,\ldots, r$, are non-negative integer numbers, $1<s_1<\ldots<s_r$.  Below we describe a representations $d_j=R_j+D_j$, $j=0,1,\ldots, r$, of the increment orders $d_j$ by stating  conditions on the fractional parts $D_j$, such that the increment sequence $$\vec y(m):=(1-B)^{R_0}\prod_{j=1}^r(1-B^{s_j})^{R_i}\vec{\xi}(m)$$ is  a stationary  fractionally integrated seasonal stochastic  sequence.
For example, in case of single  increment pattern $(1-B^{s^*})^{R^*+D^*}$, this condition is $|D^*|<1/2$.

\begin{ozn}\label{def_fract_Pryrist}
A  sequence $\chi_{\overline{s}}^{(R+D)}(\vec \xi(m))$ is called  \emph{a fractional multiple (FM) increment sequence}.
\end{ozn}

Consider the generating function of the Gegenbauer polynomial:
\[
(1-2 u B+B^2)^{-d}=\sum_{n=0}^{\infty}C_n^{(d)}(u)B^n,
\]
where
\[
C_n^{(d)}(u)=\sum_{k=0}^{[n/2]}\frac{(-1)^k(2u)^{n-2k}\Gamma(d-k+n)}{k!(n-2k)!\Gamma(d)}.
\]

The following lemma and  theorem hold true \cite{Luz_Mokl_extra_GMI}.

\begin{lema}\label{frac_incr_2}
Define the sets $\md M_j=\{\nu_{k_j}=2\pi k_j/s_j: k_j=0,1,\ldots, [s_j/2]\}$, $j=0,1,\ldots, r$, and the set $\md M=\bigcup_{j=0}^r \md M_j$. Then the increment operator $$\chi_{\overline{s}}^{(D)}(B):=(1-B)^{D_0}\prod_{j=1}^r(1-B^{s_j})^{D_j}$$
admits a representation
\begin{eqnarray*}
\chi_{\overline{s}}^{(D)}(B)
& = &\prod_{\nu \in \md M}(1-2\cos \nu B+B^2)^{\widetilde{D}_{\nu}}
\\
& = &(1-B)^{D_0+D_1+\ldots+D_r}(1+B)^{D_{\pi}}\prod_{\nu \in \md M\setminus\{0,\pi\}}(1-2\cos \nu B+B^2)^{D_{\nu}}
\\
& = &\ld(\sum_{m=0}^{\infty}G^+_{k^*}(m)B^m\rd)^{-1}=\sum_{m=0}^{\infty}G^-_{k^*}(m)B^m,
\end{eqnarray*}
where
\begin{eqnarray}
\label{Gegenbauer_GI+}
G^+_{k^*}(m)& = &\sum_{0\leq n_1,\ldots,n_{k^*}\leq m, n_1+\ldots+n_{k^*}=m}\prod_{\nu \in \md M}C_{n_{\nu}}^{(\widetilde{D}_{\nu})}(\cos\nu),
\\
\label{Gegenbauer_GI-}
G^-_{k^*}(m)& = &\sum_{0\leq n_1,\ldots,n_{k^*}\leq m, n_1+\ldots+n_{k^*}=m}\prod_{\nu \in \md M}C_{n_{\nu}}^{(-\widetilde{D}_{\nu})}(\cos\nu).
\end{eqnarray}
 $k^*=|\md M|$,  $D_{\nu}=\sum_{j=0}^rD_j \mr I \{\nu\in \md M_j\}$, $\widetilde{D}_{\nu}=D_{\nu}$ for $\nu \in \md M\setminus\{0,\pi\}$, $\widetilde{D}_{\nu}=D_{\nu}/2$ for $\nu=0$ and $\nu=\pi$.
\end{lema}

\begin{thm}\label{thm_frac}
Assume that for a stochastic vector-valued sequence $\vec \xi(m)$ and fractional differencing orders $d_j=R_j+D_j$, $j=0,1,\ldots, r$, the FM increment sequence $\chi_{\overline{1},\overline{s}}^{(R+D)}(\vec \xi(m))$ generated by increment operator (\ref{FM_increment})  is a stationary sequence with a bounded from zero and infinity spectral density $\widetilde{f}_{\overline{1}}(\lambda)$. Then for the non-negative integer numbers $R_j$, $j=0,1,\ldots, r$, the GM increment sequence $\chi_{\overline{1},\overline{s}}^{(R)}(\vec \xi(m))$    is stationary if $-1/2< D_{\nu}<1/2$ for all $\nu\in \md M$, where $D_{\nu}$ are defined by real numbers $D_j$, $j=0,1,\ldots, r$,
in Lemma \ref{frac_incr_2}, and it is long memory if $0< D_{\nu}<1/2$ for at least one $\nu\in \md M$, and invertible if $-1/2< D_{\nu}<0$. The spectral density $f(\lambda)$ of the stationary GM increment sequence $\chi_{\overline{1},\overline{s}}^{(R)}(\vec \xi(m))$ admits a representation
\[
 f(\lambda)=|\beta^{(R)}(i\lambda)|^2 \ld|\chi_{\overline{1}}^{(R)}(e^{-i\lambda})\rd|^{-2}\ld|\chi_{\overline{1}}^{(D)}(e^{-i\lambda})\rd|^{-2} \widetilde{f}_{\overline{1}}(\lambda)=:\ld|\chi_{\overline{1}}^{(D)}(e^{-i\lambda})\rd|^{-2} \widetilde{f} (\lambda),
  \]
  where
  \begin{eqnarray*}
  \ld|\chi_{\overline{1}}^{(D)}(e^{-i\lambda})\rd|^{-2}& = &\ld|\sum_{m=0}^{\infty}G^+_{k^*}(m) e^{-i\lambda m}\rd|^2=\ld|\sum_{m=0}^{\infty}G^-_{k^*}(m) e^{-i\lambda m}\rd|^{-2}
\\
& = &\prod_{\nu \in \md M}\ld|(e^{-i\nu}-e^{i\lambda})(e^{i\nu}-e^{i\lambda})\rd|^{-2\widetilde{D}_{\nu}},
  \end{eqnarray*}
 coefficients $G^+_{k^*}(m)$, $G^-_{k^*}(m)$ are defined in (\ref{Gegenbauer_GI+}), (\ref{Gegenbauer_GI-}).
\end{thm}

The spectral density $f(\lambda)$ and the structural function $D^{(R)}_{ \overline{s}}(m,\overline{1},\overline{1})$ of a stationary GM increment sequence $\chi_{\overline{1},\overline{s}}^{(R)}(\vec \xi(m))$ exhibit the following behavior for the constant matrices $C$ and $K$:
\begin{itemize}
\item   $|\beta^{(R)}(i\lambda)|^{-2} |\chi_{\overline{1}}^{(R)}(e^{-i\lambda})|^2f(\lambda)\sim C|\nu-\lambda|^{-2\widetilde{D}_{\nu}}$  as $\lambda\to \nu$, $\nu\in \mathcal{M}$ (for properties of eigenvalues of generalized fractional process, we refer to Palma and Bondon \cite{Palma-Bondon})

\item $D^{(R)}_{ \overline{s}}(m,\overline{1},\overline{1})\sim K\sum_{\nu\in \mathcal{M}:\widetilde{D}_{\nu}>0}|m|^{2\widetilde{D}_{\nu}-1}\cos (m\nu)$,   as $m\to\infty$, see Giraitis and Leipus \cite{Giraitis}.
\end{itemize}

\begin{pry}

$1.$ For the increment operator $(1-B)^{R_0+D_0}(1-B^2)^{R_1+D_1}$, $\md M_0=\{0\}$, $\md M_1=\{0,\pi\}$, $\md M=\{0,\pi\}$, the Gegenbauer representation  is $(1-B)^{D_0+D_1}(1+B)^{D_1}$.
Stationarity conditions: $|D|=|D_0+D_1|<1/2$, $|D_{\pi}|=|D_1|<1/2$.

$2.$ For the increment operator  $(1-B^2)^{R_1+D_1}(1-B^3)^{R_2+D_2}$,  $\md M_0=\{0,\pi\}$, $\md M_1=\{0,2\pi/3\}$, $\md M=\{0,2\pi/3,\pi\}$,the Gegenbauer representation   is $(1-B)^{D_1+D_2}(1-2\cos(2\pi/3)B+B^2)^{D_2}(1+B)^{D_1}$.
Stationarity conditions: $|D|=|D_1+D_2|<1/2$, $|D_{2\pi/3}|=|D_2|<1/2$, $|D_{\pi}|=|D_1|<1/2$.

$3.$ For the increment operator  $(1-B^2)^{R_1+D_1}(1-B^4)^{R_2+D_2}$, $\md M_0=\{0,\pi\}$, $\md M_1=\{0,\pi/2,\pi\}$, $\md M=\{0,\pi/2,\pi\}$, the Gegenbauer representation   is $(1-B)^{D_1+D_2}(1+B^2)^{D_2}(1+B)^{D_1+D_2}$.
Stationarity conditions: $|D|=|D_{\pi}|=|D_1+D_2|<1/2$, $|D_{\pi/2}|=|D_2|<1/2$.
\end{pry}

\section{Hilbert space projection method of interpolation}\label{classical_filtering}

\subsection{Interpolation of vector-valued stochastic sequences with stationary GM increments} \label{classical_filtering_vector}

Consider a vector-valued stochastic sequence with stationary GM increments $\vec{\xi}(m)$ constructed from transformation \eqref{PerehidXi} and a vector-valued stationary stochastic sequence $\vec\eta(m)$ uncorrelated with the sequence $\vec\eta(m)$.
Let the stationary GM increment sequence $\chi_{\overline{\mu},\overline{s}}^{(d)}(\vec{\xi}(m))=\{\chi_{\overline{\mu},\overline{s}}^{(d)}(\xi_p(m))\}_{p=1}^{T}$
and the stationary stochastic sequence $\vec\eta(m)$ have the  absolutely continuous spectral function $F(\lambda)$ and $G(\lambda)$ with the spectral densities $f(\lambda)=\{f_{ij}(\lambda)\}_{i,j=1}^{T}$ and $g(\lambda)=\{g_{ij}(\lambda)\}_{i,j=1}^{T}$ respectively.

Without loss of generality we will assume that the mean values of the increment sequences $ \mathsf{E}\chi_{\overline{\mu},\overline{s}}^{(d)}(\vec{\xi}(m))=0$, $\mathsf{E}\vec\eta(m)=0$ and $\overline{\mu}>\overline{0}$.

\textbf{Interpolation problem.} Consider the problem of mean square optimal linear estimation of the functional
\begin{equation}
A_{N}\vec{\xi}=\sum_{k=0}^{N}(\vec{a}(k))^{\top}\vec{\xi}(k)
\end{equation}
which depends on the unobserved values of the stochastic sequence $\vec{\xi}(k)=\{\xi_{p}(k)\}_{p=1}^{T}$ with stationary GM
increments.
Estimates are based on observations of the sequence $\vec\zeta(k)=\vec\xi(k)+\vec\eta(k)$  at points of the set
$\mr Z\setminus\{0,1,2\ldots,N\}$.

Assume that spectral densities $f(\lambda)$ and $g(\lambda)$ satisfy the minimality condition
\be
 \ip \text{Tr}\left[ \frac{|\beta^{(d)}(i\lambda)|^2}{|\chi_{\overline{\mu}}^{(d)}(e^{-i\lambda})|^2}\ld(f(\lambda)+|\beta^{(d)}(i\lambda)|^2 g(\lambda)\rd)^{-1}\right]
 d\lambda<\infty.
\label{umova11_f_st.n_d}
\ee

The latter condition  is the necessary and sufficient one under which the mean square errors of estimates of functional  $A\vec\xi$ is not equal to $0$.

We apply the classical Hilbert space estimation technique proposed by Kolmogorov \cite{Kolmogorov} which can be described as a $3$-stage procedure:
(i) define a target element (to be estimated) of the space $H=L_2(\Omega, \mathcal{F},\mt P)$ of random variables $\gamma$ which have zero mean values and finite variances, $\mt E\gamma=0$, $\mt E|\gamma|^2<\infty$, endowed with the inner product $\langle \gamma_1,\gamma_2\rangle={\mt E}{\gamma_1\overline{\gamma_2}}$,
(ii) define a subspace of $H$ generated by observations,
(iii) find an estimate of the target element as an orthogonal projection on the defined subspace.

\textbf{Stage i}. The  functional $A_{N}\vec{\xi}$ does not  belong to the space $H$.
With the help of the following lemma we  describe representations of the functional   as a sum of a functional with finite second moments  belonging to $H$ and a functional depending on observed values of the sequence $\vec\zeta(k)$ (``initial values'') (for more details see \cite{Luz_Mokl_filt,Luz_Mokl_book}).

\begin{lema}\label{lema predst A}
The functional $A_N\vec\zeta$ admits the representation
\be \label{zobrazh A_N_i_st.n_d}
    A_N\vec\xi=A_N\vec\zeta-A_N\vec\eta=H_N\vec\xi-V_N\vec\zeta,
\ee
where
\[
    H_N\vec\xi:=B_N\chi\vec\zeta-A_N\vec\eta,\]
    \[
    A_{N}\vec{\zeta}=\sum_{k=0}^{N}(\vec{a}(k))^{\top}\vec{\zeta}(k),\quad
    A_{N}\vec{\eta}=\sum_{k=0}^{N}(\vec{a}(k))^{\top}\vec{\eta}(k),\]
\[
B_{N}\chi\vec{\zeta}=\sum_{k=0}^{N}(\vec{b}_{N}(k))^{\top}\chi_{\overline{\mu},\overline{s}}^{(d)}(\vec{\zeta}(k)),
\quad
V_{N}\vec{\zeta}=\sum_{k=-n(\gamma)}^{-1} (\vec{v}_{N}(k))^{\top}\vec{\zeta}(k),
\]
the coefficients
$\vec{b}_{N}(k)=\{b_{N,p}(k)\}_{p=1}^{T}, k=0,1,\dots,N$ and
 $\vec{v}_{N}(k)=\{ v_{N,p}(k)\}_{p=1}^{T},$ $ k=-1,-2,\dots,-n(\gamma)$
are calculated by the formulas
\begin{eqnarray}\label{koefv_N_diskr}
 \vec{v}_N(k)& = &\sum_{l=0}^{N \wedge k+n(\gamma)} \mt{diag}_T(e_{\nu}(l-k))\vec{b}_N(l), \, k=-1,-2,\dots,-n(\gamma),
\\
 \label{determ_b_N}
\vec{b}_N(k)& = &\sum_{m=k}^N\mt{diag}_T(d_{\overline{\mu}}(m-k))\vec{a}(m)=(D^{\overline{\mu}}_{N}{\me a}_{N})_k, \, k=0,1,\dots,N,
\end{eqnarray}
  $D^{\overline{\mu}}_{N}$ is the  linear transformation  determined by a matrix with the entries $(D^{\overline{\mu}}_{N})(k,j)=\mt{diag}_T(d_{\overline{\mu}}(j-k))$ if
$0\leq k\leq j\leq N$, and $(D^{\overline{\mu}}_{N})(k,j)=0$ if $0\leq j<k\leq N$, $\mt{diag}_T(x)$ denotes a $T\times T$ diagonal matrix with the entry $x$ on its diagonal, $\me a_N=((\vec{a}(0))^{\top},(\vec{a}(1))^{\top}, \ldots,(\vec{a}(N))^{\top})^{\top}$, coefficients $\{d_{\overline{\mu}}(k):k\geq0\}$ are determined by the relationship
\[
 \sum_{k=0}^{\infty}d_{\overline{\mu}}(k)x^k
=\prod_{i=1}^r\ld(\sum_{j_i=0}^{\infty}x^{\mu_is_ij_i}\rd)^{d_i}.\]

\end{lema}

The functional $H_N\vec\xi$ from representation (\ref{zobrazh A_N_i_st.n_d}) has finite variance and the functional  $V_N\vec\zeta$ depends on the known observations of the stochastic sequence $\vec\zeta(k)$ at points $k=-n(\gamma),-n(\gamma)+1,\ldots,-1$. Therefore, estimates $\widehat{A}_N\vec\xi$ and $\widehat{H}_N\vec\xi$ of the functionals $A_N\vec\xi$ and $H_N\vec\xi$ and the mean-square errors $\Delta(f,g;\widehat{A}_N\vec\xi)=\mt E |A_N\vec\xi-\widehat{A}_N\vec\xi|^2$ and $\Delta(f,g;\widehat{H}_N\vec\xi)=\mt E
|H_N\vec\xi-\widehat{H}_N\vec\xi|^2$ of the estimates $\widehat{A}_N\vec\xi$ and $\widehat{H}_N\vec\xi$ satisfy the following relations
\be\label{mainformula}
    \widehat{A}_N\vec\xi=\widehat{H}_N\vec\xi-V_N\vec\zeta,\ee
\be\label{mainformula2}
    \Delta(f,g;\widehat{A}_N\vec\xi)
    =\mt E |A_N\vec\xi-\widehat{A}_N\vec\xi|^2=
  \mt E|H_N\vec\xi-\widehat{H}_N\vec\xi|^2=\Delta(f,g;\widehat{H}_N\vec\xi).
    \ee
Therefore, the estimation problem for the functional $A_N\vec\xi$ is equivalent to the estimation problem for the functional $H_N\vec\xi$.
This problem can be solved by applying the Hilbert space projection method proposed by Kolmogorov \cite{Kolmogorov}.

With the help of the spectral representations of stochastic sequences involved we can write the following spectral representation of the  functional
\begin{equation*}
H_N\vec\xi=
\int_{-\pi}^{\pi}
\left(\vec{B}_{\overline{\mu},N}(e^{i\lambda})\right)^{\top}
\frac{\chi_{\overline{\mu}}^{(d)}(e^{-i\lambda})}{\beta^{(d)}(i\lambda)}
d\vec{Z}_{\xi^{(d)}+\eta^{(d)}}(\lambda)-
\int_{-\pi}^{\pi}\left(\vec{A}_N(e^{i\lambda})\right)^{\top}d\vec{Z}_{\eta}(\lambda),
\end{equation*}
where
\begin{equation*}
\vec{B}_{\overline{\mu},N}(e^{i\lambda})=\sum_{k=0}^{N}\vec{b}_{\overline{\mu},N}(k)e^{i\lambda k}
=\sum_{k=0}^{N}(D^{\overline{\mu}}_N\me a_N)_ke^{i\lambda k},\quad
 \vec{A}_N(e^{i\lambda})=\sum_{k=0}^{N}\vec{a}(k)e^{i\lambda k}.
\end{equation*}

\emph{Stage (ii).} Introduce the following notations. Denote by
$H^{0-}(\xi^{(d)}_{\overline{\mu},\overline{s}}+\eta^{(d)}_{\overline{\mu},\overline{s}})$ the closed linear subspace generated by values
$\{\chi_{\overline{\mu},\overline{s}}^{(d)}(\vec{\xi}(k))+\chi_{\overline{\mu},\overline{s}}^{(d)}(\vec{\eta}(k)):k=-1,-2,-3,\dots\}$, $\overline{\mu}>\vec 0$
of the observed GM increments
in the Hilbert space $H=L_2(\Omega,\mathcal{F},\mt P)$ of random variables $\gamma$ with zero mean value, $\mt E\gamma=0$, finite variance, $\mt E|\gamma|^2<\infty$, and the inner product $(\gamma_1;\gamma_2)=\mt E\gamma_1\overline{\gamma_2}$.

Denote by
$H^{N+}(\xi^{(d)}_{-\overline{\mu},\overline{s}}+\eta^{(d)}_{-\overline{\mu},\overline{s}})$ the closed linear subspace generated by values of the sequence
$\{\chi_{- \overline{\mu},\overline{s}}^{(d)}(\vec{\xi}(k))+\chi_{-\overline{\mu},\overline{s}}^{(d)}(\vec{\eta}(k)):k\geq N\}$,
of the observed GM increments
in the Hilbert space $H=L_2(\Omega,\mathcal{F},\mt P)$.

Denote by $L_2^{0-}(f(\lambda)+|{\beta^{(d)}(i\lambda)}|^2 g(\lambda))$ and $L_2^{N+}(f(\lambda)+ |{\beta^{(d)}(i\lambda)}|^2 g(\lambda))$ the closed linear subspaces of the Hilbert space
$L_2(f(\lambda)+|{\beta^{(d)}(i\lambda)}|^2 g(\lambda))$  of vector-valued functions with the inner product $\langle g_1;g_2\rangle=\ip (g_1(\lambda))^{\top}(f(\lambda)+|{\beta^{(d)}(i\lambda)}|^2 g(\lambda))\overline{g_2(\lambda)}d\lambda$ which is generated by the functions
\[
 e^{i\lambda k}\chi_{\overline{\mu}}^{(d)}(e^{-i\lambda})\frac{1}{\beta^{(d)}(i\lambda)}\vec\delta_l,\quad \vec\delta_l=\{\delta_{lp}\}_{p=1}^T,\,
 l=1,\dots,T;\,\text{for}\,\, k \leq -1,\,\,\text{and} \,\,k\geq N+1,
\]
respectively, where $\delta_{lp}$ are Kronecker symbols.

Then the relation
\[
 \chi_{\overline{\mu},\overline{s}}^{(d)}(\vec{\xi}(k))+\chi_{\overline{\mu},\overline{s}}^{(d)}(\vec{\eta}(k))=\ip e^{i\lambda
 k}\chi_{\overline{\mu}}^{(d)}(e^{-i\lambda})\dfrac{1}{\beta^{(d)}(i\lambda)}dZ_{\xi^{(d)}+\eta^{(d)}}(\lambda)\]
yields a  one to one correspondence between elements
\[
e^{i\lambda k}\chi_{\overline{\mu}}^{(d)}(e^{-i\lambda})\dfrac{1}{\beta^{(d)}(i\lambda)}\vec \delta_l\]
from the space

$$L_2^{0-}(f(\lambda)+|{\beta^{(d)}(i\lambda)}|^2 g(\lambda))  \oplus L_2^{N+}(f(\lambda)+|{\beta^{(d)}(i\lambda)}|^2 g(\lambda))  $$
and elements
$\chi_{\overline{\mu},\overline{s}}^{(d)}(\vec{\xi}(k))+\chi_{\overline{\mu},\overline{s}}^{(d)}(\vec{\eta}(k))$
from the space
\[
H^{0-}(\xi^{(d)}_{\overline{\mu},\overline{s}}+\eta^{(d)}_{\overline{\mu},\overline{s}})
\oplus
H^{N+}(\xi^{(d)}_{-\overline{\mu},\overline{s}}+\eta^{(d)}_{-\overline{\mu},\overline{s}})
=
H^{0-}(\xi^{(d)}_{\overline{\mu},\overline{s}}+\eta^{(d)}_{\overline{\mu},\overline{s}}) \oplus
 H^{(N+ n(\gamma))+}(\xi^{(d)}_{\overline{\mu},\overline{s}}+\eta^{(d)}_{\overline{\mu},\overline{s}}).
\]

Relation (\ref{mainformula}) implies that every linear estimate $\widehat{A}\vec\xi$ of the functional $A\vec\xi$
can be represented in the form
\begin{equation}
 \label{otsinka A_e_d}
 \widehat{A}_N\vec\xi=\ip
(\vec{h}_{\overline{\mu},N}(\lambda))^{\top}d\vec{Z}_{\xi^{(d)}+\eta^{(d)}}(\lambda)-
\sum_{k=-\mu n}^{-1}(\vec v_{N}(k))^{\top}(\vec\xi(k)+\vec\eta(k)),
\end{equation}
 where
$\vec{h}_{\overline{\mu},N}(\lambda)=\{h_{p}(\lambda)\}_{p=1}^{T}$ is the spectral characteristic of the optimal estimate $\widehat{H}_N\vec\xi$.

\emph{Stage (iii).}
At this stage we find the mean square optimal estimate
$\widehat{H}_N\vec\xi$ as a projection of the element $H_N\vec\xi$ on the
subspace
$H^{0-}(\xi^{(d)}_{\overline{\mu},\overline{s}}+\eta^{(d)}_{\overline{\mu},\overline{s}}) \oplus H^{(N+ n(\gamma))+}(\xi^{(d)}_{\overline{\mu},\overline{s}}+\eta^{(d)}_{\overline{\mu},\overline{s}})$.  This projection is
determined by two conditions:

1) $ \widehat{H}_N\vec\xi\in H^{0-}(\xi^{(d)}_{\overline{\mu},\overline{s}}+\eta^{(d)}_{\overline{\mu},\overline{s}}) \oplus H^{(N+ n(\gamma))+}(\xi^{(d)}_{\overline{\mu},\overline{s}}+\eta^{(d)}_{\overline{\mu},\overline{s}})$;

2) $(H_N\vec\xi-\widehat{H}_N\vec\xi)
\perp
H^{0-}(\xi^{(d)}_{\overline{\mu},\overline{s}}+\eta^{(d)}_{\overline{\mu},\overline{s}}) \oplus H^{(N+ n(\gamma))+}(\xi^{(d)}_{\overline{\mu},\overline{s}}+\eta^{(d)}_{\overline{\mu},\overline{s}})$.

The second condition implies the following relation which holds true for all $k\leq-1$ and $k\geq N+n(\gamma) +1$
\begin{multline*}
\int_{-\pi}^{\pi}
\bigg[\bigg(
\left(\vec{B}_{\overline{\mu},N}(e^{i\lambda})\right)^{\top}
\frac{\chi_{\overline{\mu}}^{(d)}(e^{-i\lambda})}{\beta^{(d)}(i\lambda)}
-
\vec{h}_{\overline{\mu},N}(\lambda)\bigg)^{\top}
(f(\lambda)+|{\beta^{(d)}(i\lambda)}|^2 g(\lambda))
-
\\
-
(\vec{A}_N(e^{i\lambda})^{\top}g(\lambda){\overline{\beta^{(d)}(i\lambda)}}\bigg]
\frac{\chi_{\overline{\mu}}^{(d)}(e^{-i\lambda})}
{\overline{\beta^{(d)}(i\lambda)}}
e^{-i\lambda k}d\lambda=0.
 \end{multline*}

\noindent This relation allows us to derive the spectral characteristic
$\vec{h}_{\overline{\mu},N}(\lambda)$ of the estimate $\widehat{H}_N\vec\xi$ which can be represented in the form
\begin{multline} \label{spectr A}
(\vec{h}_{\overline{\mu},N}(\lambda))^{\top}=
\left(\vec{B}_{\overline{\mu},N}(e^{i\lambda})\right)^{\top}
\frac{\chi_{\overline{\mu}}^{(d)}(e^{-i\lambda})}{\beta^{(d)}(i\lambda)}
-
\\
-
(\vec{A}_N(e^{i\lambda})^{\top}g(\lambda){\overline{\beta^{(d)}(i\lambda)}}
(f(\lambda)+|{\beta^{(d)}(i\lambda)}|^2 g(\lambda))
^{-1}
-
\\
-(\vec{C}_{\overline{\mu},N}(e^{i\lambda}))^{\top}
\frac{\overline{\beta^{(d)}(i\lambda)}}{\chi_{\overline{\mu}}^{(d)}(e^{-i\lambda})}
(f(\lambda)+|{\beta^{(d)}(i\lambda)}|^2 g(\lambda))^{-1},
\end{multline}
\[
\vec{C}_{\overline{\mu},N}(e^{i \lambda})=\sum_{k=0}^{N+n(\gamma)}\vec{c}_{\overline{\mu},N}(k)e^{ik\lambda},\]
 $\vec{c}_{\overline{\mu},N}(k)=\{c_{\overline{\mu},N,p}(k)\}_{p=1}^T, k=0,1,\dots,N+n(\gamma),$  are unknown coefficients to be found.

It follows from condition 1) that the following equations should be satisfied  for $0\leq j\leq N+n(\gamma)$
\begin{multline} \label{eq_C}
\int_{-\pi}^{\pi} \biggl[(\vec{B}_{\overline{\mu},N}(e^{i\lambda}))^{\top}-
(\vec{A}_N(e^{i\lambda})^{\top}g(\lambda)
\frac{|{\beta^{(d)}(i\lambda)}|^2}{\chi_{\overline{\mu}}^{(d)}(e^{-i\lambda})}
(f(\lambda)+|{\beta^{(d)}(i\lambda)}|^2 g(\lambda))^{-1}
-
\\
-
(\vec{C}_{\overline{\mu},N}(e^{i\lambda}))^{\top}
\frac{|{\beta^{(d)}(i\lambda)}|^2}{|\chi_{\overline{\mu}}^{(d)}(e^{-i\lambda})|^2}
(f(\lambda)+|{\beta^{(d)}(i\lambda)}|^2 g(\lambda))^{-1}
\biggr]e^{-ij\lambda}d\lambda=0.
\end{multline}

 Define for $0\leq k,j\leq N+n(\gamma)$ the Fourier coefficients of the corresponding functions
\[
T^{\overline{\mu}}_{k,j}=\frac{1}{2\pi}\int_{-\pi}^{\pi}
e^{i\lambda(j-k)}
 \frac{|\beta^{(d)}(i\lambda)|^2}{\chi_{\overline{\mu}}^{(d)}(e^{-i\lambda})}g(\lambda)
(f(\lambda)+|{\beta^{(d)}(i\lambda)}|^2 g(\lambda))^{-1}
d\lambda;
\]
\[
P_{k,j}^{\overline{\mu}}=\frac{1}{2\pi}\int_{-\pi}^{\pi} e^{i\lambda (j-k)}
 \frac{|\beta^{(d)}(i\lambda)|^2}{|\chi_{\overline{\mu}}^{(d)}(e^{-i\lambda})|^2}
(f(\lambda)+|{\beta^{(d)}(i\lambda)}|^2 g(\lambda))^{-1}
d\lambda;
\]
 \[
 Q_{k,j}=\frac{1}{2\pi}\int_{-\pi}^{\pi}
e^{i\lambda(j-k)}f(\lambda)g(\lambda)(f(\lambda)+|{\beta^{(d)}(i\lambda)}|^2 g(\lambda))^{-1}
d\lambda.
\]

  \noindent Making use of the defined Fourier coefficients, relation \eqref{eq_C} can be presented as a system of $N+n(\gamma)+1$ linear equations
  determining the unknown coefficients ${\vec c}_{{\overline{\mu}},N}(k)$, $0\leq k\leq N+\mu n$.
\begin{equation} \label{linear equations1}
    \vec{b}_{\overline{\mu},N}(j)-\sum_{ m=0}^{N+n(\gamma)}T^{\overline{\mu}}_{j,m}\vec{a}_{\overline{\mu},N}(m)
    =\sum_{k=0}^{N+n(\gamma)}P_{j,k}^{\overline{\mu}}\vec{c}_{\overline{\mu},N}(k),\,\, 0\leq j\leq N,
    \end{equation}
\begin{equation}  \label{linear equations2}
    -\sum_{ m=0}^{N+n(\gamma)}T^{\overline{\mu}}_{j,m}\vec{a}_{\overline{\mu},N}(m)
    =\sum_{k=0}^{N+n(\gamma)}P_{j,k}^{\overline{\mu}}\vec{c}_{\overline{\mu},N}(k),\,\, N+1\leq j\leq N+n(\gamma),
    \end{equation}
where coefficients $\{\vec{a}_{\overline{\mu},N}(m):0\leq m\leq N+n(\gamma)\}$ are calculated by the formula
\[
   \vec   a_{\overline{\mu},N}(m)=\vec a_{-\overline{\mu},N}(m- n(\gamma)),\quad 0\leq m\leq N+n(\gamma),\]
\begin{equation} \label{coeff a_N_mu}
 \vec a_{-\overline{\mu},N}(m)=\sum_{l=\max\ld\{m,0\rd\}}^{\min\{m+n(\gamma),N\}}e_{\gamma}(l-m)\vec a(l),\quad  -n(\gamma)\leq m\leq N.
 \end{equation}

  \noindent Denote by $[D_N^{\overline{\mu}}\me a_N]_{+n(\gamma)}$ a vector of dimension  $(N+n(\gamma) +1)T$ which is constructed by adding
$(n(\gamma))T$ zeros to the vector $D_N^{\overline{\mu}}\me a_N$ of dimension  $(N+1)T$.
Making use of this definition the system (\ref{linear equations1}) -- (\ref{linear equations2}) can be represented in the matrix form:
\[[D_N^{\overline{\mu}}\me a_N]_{+n(\gamma)}-\me T^{\overline{\mu}}_N\me a^{\overline{\mu}}_N=\me P^{\overline{\mu}}_N\me c^{\overline{\mu}}_N,\]
where
\[
    \me a^{\overline{\mu}}_N=((\vec{a}_{\overline{\mu},N}(0))^{\top},(\vec{a}_{{\overline{\mu}},N}(1))^{\top},(\vec{a}_{{\overline{\mu}},N}(2))^{\top},\ldots,(\vec{a}_{{\overline{\mu}},N}(N+n(\gamma)))^{\top})^{\top}
    \]
 \[
    \me c^{{\overline{\mu}}}_N=((\vec{c}_{{\overline{\mu}},N}(0))^{\top},(\vec{c}_{{\overline{\mu}},N}(1))^{\top},(\vec{c}_{{\overline{\mu}},N}(2))^{\top},\ldots,(\vec{c}_{{\overline{\mu}},N}(N+n(\gamma)))^{\top})^{\top}
    \]
are  vectors of  dimension $(N+n(\gamma)+1)T$; $\me P^{{\overline{\mu}}}_N$ and $\me T^{{\overline{\mu}}}_N$ are matrices of dimension $(N+n(\gamma)+1)T\times(N+n(\gamma)+1)T$ with $T\times T$ matrix elements $(\me P^{{\overline{\mu}}}_N)_{j,k}=P_{j,k}^{{\overline{\mu}}}$ and $(\me T^{{\overline{\mu}}}_N)_{j, k} =T^{{\overline{\mu}}}_{j,k}$, $0\leq j,k\leq N+n(\gamma)$.

Thus, the coefficients $\vec{c}_{{\overline{\mu}},N}(k)$, $0\leq k\leq N+n(\gamma)$, are determined by the formula $ 0\leq k\leq N+n(\gamma)$
\[
  \vec{c}_{{\overline{\mu}},N}(k)=\ld((\me P^{{\overline{\mu}}}_N)^{-1}[D_N^{{\overline{\mu}}}\me a_N]_{+n(\gamma)}-(\me P^{{\overline{\mu}}}_N)^{-1}\me T^{{\overline{\mu}}}_N\me a^{{\overline{\mu}}}_{N}\rd)_k,\]
where $\ld((\me P^{{\overline{\mu}}}_N)^{-1}[D_N^{{\overline{\mu}}}\me a_N]_{+n(\gamma)}-(\me P^{{\overline{\mu}}}_N)^{-1}\me T^{{\overline{\mu}}}_N\me a^{{\overline{\mu}}}_N\rd)_k$, $0\leq k\leq N+n(\gamma)$, is the
 $k$th element of the vector  $(\me P^{{\overline{\mu}}}_N)^{-1}[D_N^{{\overline{\mu}}}\me a_N]_{+n(\gamma)}-(\me P^{{\overline{\mu}}}_N)^{-1}\me T^{{\overline{\mu}}}_N\me a^{{\overline{\mu}}}_N$.

The existence of the inverse matrix $(\me P^{\mu}_N)^{-1}$ was shown in \cite{Luz_Mokl_book} under condition (\ref{umova11_f_st.n_d}).

The spectral characteristic  $\vec{h}_{{\overline{\mu}},N}(\lambda)$ of the estimate $\widehat{H}_N\xi$ of the functional $H_N\xi$ is calculated by formula (\ref{spectr A}), where
\begin{equation} \label{spectr C}
 \vec{C}_{{\overline{\mu}},N}(e^{i \lambda})=\sum_{k=0}^{N+n(\gamma)}
    \bigg((\me P^{{\overline{\mu}}}_N)^{-1}[D_N^{{\overline{\mu}}}\me a_N]_{+n(\gamma)}- (\me P^{{\overline{\mu}}}_N)^{-1}\me T^{{\overline{\mu}}}_N\me a^{{\overline{\mu}}}_N\bigg)_k e^{i\lambda k}.
\end{equation}

The value of the  mean-square errors of the estimates $\widehat{A}_N\vec\xi$ and $\widehat{H}_N\vec\xi$ can be calculated by the formula
\[
\Delta(f,g;\widehat{A}_N\vec\xi)=\Delta(f,g;\widehat{H}_N\vec\xi)= \mt E|H_N\vec\xi-\widehat{H}_N\vec\xi|^2=
\]
\[
=
\frac{1}{2\pi}\int_{-\pi}^{\pi}
\frac{|\beta^{(d)}(i\lambda)|^2}{|\chi_{\overline{\mu}}^{(d)}(e^{-i\lambda})|^2}
\left[\chi_{\overline{\mu}}^{(d)}(e^{i\lambda})(\vec{A}_N(e^{i\lambda}))^{\top}g(\lambda) -
(\vec{C}_{\overline{\mu},N}(e^{i \lambda}))^{\top}
\right]
\times
\]
\[
\times
(f(\lambda)+|\beta^{(d)}(i\lambda)|^2g(\lambda))^{-1}\, f(\lambda)\, (f(\lambda)+|\beta^{(d)}(i\lambda)|^2 g(\lambda))^{-1}
\times
\]
\[
\times
\left[\chi_{\overline{\mu}}^{(d)}(e^{i\lambda}){\vec{A}_N(e^{-i\lambda})}g(\lambda) -
(\vec{C}_{\overline{\mu},N}(e^{-i \lambda}))
\right]
d\lambda+
\]
\[
+\frac{1}{2\pi}\int_{-\pi}^{\pi}
 \frac{1}{|\chi_{\overline{\mu}}^{(d)}(e^{-i\lambda})|^2}
\left[\chi_{\overline{\mu}}^{(d)}(e^{i\lambda})(\vec{A}_N(e^{i\lambda}))^{\top}f(\lambda) + |\beta^{(d)}(i\lambda)|^2
(\vec{C}_{\overline{\mu},N}(e^{i \lambda}))^{\top}
\right]
\times
\]
\[
\times
(f(\lambda)+|\beta^{(d)}(i\lambda)|^2 g(\lambda))^{-1}\,g(\lambda)\,(f(\lambda)+|\beta^{(d)}(i\lambda)|^2 g(\lambda))^{-1}
\times
\]
\[
\times
\left[\chi_{\overline{\mu}}^{(d)}(e^{i\lambda}){\vec{A}_N(e^{-i\lambda})}f(\lambda) + |\beta^{(d)}(i\lambda)|^2
(\vec{C}_{\overline{\mu},N}(e^{-i \lambda}))
\right]
d\lambda=
\]
 \be
= \ld\langle [D_N^{\overline{\mu}}\me a_N]_{+n(\gamma)}- \me T^{\overline{\mu}}_N\me
    a^{\overline{\mu}}_N,(\me P^{\overline{\mu}}_N)^{-1}[D_N^{\overline{\mu}}\me a_N]_{+n(\gamma)}
    -(\me P^{\overline{\mu}}_N)^{-1}\me T^{\overline{\mu}}_N\me a^{\overline{\mu}}_N\rd\rangle
    \\+\langle\me Q_N\me a_N,\me
    a_N\rangle,
    \label{pohybkaAN}
    \ee
where $\me Q_N$ is a matrix of the dimension $(N+1)T\times(N+1)T$ with the $T\times T$ matrix elements $(\me Q_N)_{j,k}=Q_{j,k}$, $0\leq j,k\leq N$

\begin{thm}\label{thm_intA}
Let $\{\vec\xi(m),m\in\mr Z\}$ be a stochastic sequence which defines
the stationary GM increment sequence $\chi_{\overline{\mu},\overline{s}}^{(d)}(\vec{\xi}(m))=\{\chi_{\overline{\mu},\overline{s}}^{(d)}(\xi_p(m))\}_{p=1}^{T}$ with the
absolutely continuous spectral function $F(\lambda)$ which has
spectral density $f(\lambda)$. Let $\{\vec\eta(m),m\in\mr Z\}$ be an
uncorrelated with the sequence $\vec\xi(m)$ stationary stochastic
sequence with an absolutely continuous spectral function
$G(\lambda)$ which has spectral density $g(\lambda)$. Let the minimality condition
(\ref{umova11_f_st.n_d}) be satisfied. The
optimal linear estimate $\widehat{A}_N\vec\xi$ of the functional
$A_N\vec\xi$ which depends on the unknown values of elements $\vec\xi(k)$, $k=0,1,2,\ldots,N$,   from
observations of the sequence $\vec\xi (m )+\vec\eta (m )$ at points of the set $Z\setminus\{0,1,2,\ldots,N\}$
is calculated by  formula (\ref{otsinka A_e_d}).
The spectral characteristic $\vec{h}_{{\overline{\mu}},N}(\lambda)$ of the optimal
estimate $\widehat{A}_N\vec\xi$ is calculated by  formulas (\ref{spectr A}), (\ref{spectr C}). The value of the mean-square error
$\Delta(f,g;\widehat{A}_N\vec\xi)$ is calculated by  formula (\ref{pohybkaAN}).
\end{thm}

\begin{nas}
The spectral characteristic $\vec{h}_{{\overline{\mu}},N}(\lambda)$ (\ref{spectr A}) admits the representation
$\vec{h}_{{\overline{\mu}},N}(\lambda)=\vec{h}_{{\overline{\mu}},N}^{1}(\lambda)-\vec{h}_{{\overline{\mu}},N}^2(\lambda)$, where
\begin{multline}
\label{spectr h1}
(\vec{h}_{{\overline{\mu},N}}^1(\lambda))^{\top}=(\vec{B}_{{\overline{\mu}},N}(e^{i\lambda}))^{\top}
\frac{\chi_{\overline{{{\mu}}}}^{(d)}(e^{-i\lambda})}{\beta^{(d)}(i\lambda)}-
\\
-\frac{\overline{\beta^{(d)}(i\lambda)}}{\chi_{\overline{{{\mu}}}}^{(d)}(e^{-i\lambda})}
\left(
\sum_{k=0}^{N+n(\gamma)}
    \ld((\me P^{{\overline{\mu}}}_N)^{-1}[D_N^{{\overline{\mu}}}\me a_N]_{{+n(\gamma)}}\rd)_k e^{i\lambda k}
\right)^{\top}
(f(\lambda)+|\beta^{(d)}(i\lambda)|^2 g(\lambda))^{-1},
 \end{multline}
 \begin{multline} \label{spectr h2}
(\vec{h}_{{\overline{\mu}},N}^2(\lambda))^{\top}=
(\vec{A}_N(e^{i\lambda }))^{\top}
 {{\overline{\beta^{(d)}(i\lambda)}} g(\lambda)}(f(\lambda)+|\beta^{(d)}(i\lambda)|^2 g(\lambda))^{-1} -
 \\
-\frac{\overline{\beta^{(d)}(i\lambda)}}{\chi_{\overline{{{\mu}}}}^{(d)}(e^{-i\lambda})}
\left(
\sum_{k=0}^{N+n(\gamma)}
    \ld((\me P^{{\overline{\mu}}}_N)^{-1}\me T^{{\overline{\mu}}}_N\me a^{{\overline{\mu}}}_N\rd)_k e^{i\lambda k}
\right)^{\top}
(f(\lambda)+|\beta^{(d)}(i\lambda)|^2 g(\lambda))^{-1}.
\end{multline}
Here $\vec{h}_{\overline{\mu},N}^1(\lambda)$ and $\vec{h}_{\overline{\mu},N}^2(\lambda)$ are spectral characteristics of the optimal estimates
$\widehat{B}_N\chi\vec\zeta$ and $\widehat{A}_N\vec\eta$ of the functionals $B_N\chi\vec\zeta$ and $A_N\vec\eta$ respectively based on observations
$\vec\xi(k)+\vec\eta(k)$ at points of the set $Z\setminus\{0,1,2,\ldots,N\}$.
\end{nas}

\begin{zau}
The interpolation problem for stochastic sequences with fractional multiple (FM) increments can be solved with the help of results described in Theorem \ref{thm_intA} under the conditions of Theorem \ref{thm_frac} on the increment  orders $d_i$.
\end{zau}

\subsection{Interpolation of stochastic sequences with periodically  stationary GM increments} \label{classical_filtering_vector}

Consider the problem of mean square optimal linear estimation of the functional
\begin{equation}
A_{M}{\vartheta}=\sum_{k=0}^{N}{a}^{(\vartheta)}(k)\vartheta(k)
\end{equation}
which depend on unobserved values of the stochastic sequence ${\vartheta}(m)$ with with periodically stationary GM increments.
Estimates are based on observations of the sequence $\zeta(m)=\vartheta(m)+\eta(m)$ at points of the set $Z\setminus\{0,1,2,\ldots,N\}$,
 where the periodically stationary noise sequence $ \eta(m)$ is uncorrelated with ${\vartheta}(m)$.

The functional $A_M{\vartheta}$ can be represented in the form
\begin{eqnarray}
\nonumber
A_M{\vartheta}& = &\sum_{k=0}^{M}{a}^{(\vartheta)}(k)\vartheta(k)=\sum_{m=0}^{N}\sum_{p=1}^{T}
{a}^{(\vartheta)}(mT+p-1)\vartheta(mT+p-1)
\\\nonumber
& = &\sum_{m=0}^{N}\sum_{p=1}^{T}a_p(m)\xi_p(m)=\sum_{m=0}^{N}(\vec{a}(m))^{\top}\vec{\xi}(m)=A_N\vec{\xi},
\end{eqnarray}
where $N=[\frac{M}{T}]$, the sequence $\vec{\xi}(m) $ is determined by the formula
\be \label{zeta}
\vec{\xi}(m)=({\xi}_1(m),{\xi}_2(m),\dots,{\xi}_T(m))^{\top},\,
 {\xi}_p(m)=\vartheta(mT+p-1);\,p=1,2,\dots,T;
\ee
\begin{eqnarray}
\nonumber
(\vec{a}(m))^{\top}& = &({a}_1(m),{a}_2(m),\dots,{a}_T(m))^{\top},
\\\nonumber
 {a}_p(m)& = &a^{\vartheta}(mT+p-1);\,0\leq m\leq N; 1\leq p\leq T;\,mT+p-1\leq M;
\\  {a}_p(N)& = &0;\quad
M+1\leq NT+p-1\leq (N+1)T-1;1\leq p\leq T. \label{aNzeta}
\end{eqnarray}

Making use of the introduced notations and statements of Theorem \ref{thm_intA} we can claim that the following theorem holds true.

\begin{thm}
\label{thm_est_A_Nzeta}
Let a stochastic sequence ${\vartheta}(k)$ with periodically stationary GM increments generate by formula \eqref{zeta}
 a vector-valued stochastic sequence $\vec{\xi}(m) $ which determine a
stationary  GM  increment sequence
$\chi_{\overline{\mu},\overline{s}}^{(d)}(\vec{\xi}(m))=\{\chi_{\overline{\mu},\overline{s}}^{(d)}(\xi_p(m))\}_{p=1}^{T}$
with the spectral density matrix $f(\lambda)=\{f_{ij}(\lambda)\}_{i,j=1}^{T}$.
Let $\{\vec\eta(m),m\in\mr Z\}$, $\vec{\eta}(m)=({\eta}_1(m),{\eta}_2(m),\dots,{\eta}_T(m))^{\top},\,
 {\eta}_p(m)=\eta(mT+p-1);\,p=1,2,\dots,T, $
be uncorrelated with the sequence $\vec\xi(m)$ stationary stochastic
sequence with an absolutely continuous spectral function
$G(\lambda)$ which has spectral density $g(\lambda)$. Let the minimality condition
(\ref{umova11_f_st.n_d}) be satisfied.
Let coefficients $\vec {a}(k), k\geqslant 0$ be determined by formula \eqref{aNzeta}.
The optimal linear estimate $\widehat{A}_M\zeta$ of the functional $A_M\zeta=A_N\vec{\xi}$ based on observations of the sequence
$\zeta(m)=\vartheta(m)+\eta(m)$ at points of the set $Z\setminus\{0,1,2,\ldots,N\}$ is calculated by
 formula (\ref{otsinka A_e_d}).
The spectral characteristic $\vec{h}_{\overline{\mu},N}(\lambda)=\{h_{\overline{\mu},N,p}(\lambda)\}_{p=1}^{T}$
and the value of the mean square error $\Delta(f;\widehat{A}_M\zeta)$
are calculated by formulas  (\ref{spectr A}), (\ref{spectr C}) and \eqref{pohybkaAN} respectively.
\end{thm}

\section{Minimax (robust) method of interpolation}\label{minimax_filtering}

The values of the mean square errors and the spectral characteristics of the optimal estimate
of the functional ${A}_N\vec\xi$
depending on the unobserved values of a stochastic sequence $\vec{\xi}(m)$ which determine a stationary  GM increments sequence
$\chi_{\overline{\mu},\overline{s}}^{(d)}(\vec{\xi}(m))$
with the spectral density matrix $f(\lambda)$
based on observations of the sequence
$\vec\xi(m)+\vec\eta(m)$ at points $Z\setminus\{0,1,2,\ldots,N\}$ can be calculated by formulas
(\ref{spectr A}), (\ref{spectr C}), \eqref{pohybkaAN}
respectively,
in the case where the spectral density matrices
$f(\lambda)$ and $g(\lambda)$ of the target sequence  and the noise are exactly known.

In practical cases, however,  complete information about the spectral density matrices is not available.
If in such cases a set $\md D=\md D_f\times\md D_g$ of admissible spectral densities is defined,
the minimax-robust approach to
estimation of linear functionals depending on unobserved values of stochastic sequences with stationary increments may be applied.

This method consists in finding an estimate that minimizes
the maximal values of the mean square errors for all spectral densities
from a given class $\md D=\md D_f\times\md D_g$ of admissible spectral densities
simultaneously. This method will be applied in the case of concrete classes of spectral densities.

To formalize this approach we recall the following definitions \cite{Moklyachuk}.

\begin{ozn}
For a given class of spectral densities $\mathcal{D}=\md
D_f\times\md D_g$ the spectral densities
$f_0(\lambda)\in\mathcal{D}_f$, $g_0(\lambda)\in\md D_g$ are called
least favorable in the class $\mathcal{D}$ for the optimal linear
estimation of the functional $A_N\vec \xi$  if the following relation holds
true:
\[\Delta(f_0,g_0)=\Delta(h(f_0,g_0);f_0,g_0)=
\max_{(f,g)\in\mathcal{D}_f\times\md
D_g}\Delta(h(f,g);f,g).\]
\end{ozn}

\begin{ozn}
For a given class of spectral densities $\mathcal{D}=\md
D_f\times\md D_g$ the spectral characteristic $h^0(\lambda)$ of
the optimal linear estimate of the functional $A_N\vec \xi$ is called
minimax-robust if there are satisfied the conditions
\[h^0(\lambda)\in H_{\mathcal{D}}=\bigcap_{(f,g)\in\mathcal{D}_f\times\md D_g}L_2^{0-}
(f(\lambda)+|\beta^{(d)}(i\lambda)|^2 g(\lambda)),
\]
\[\min_{h\in H_{\mathcal{D}}}\max_{(f,g)\in \mathcal{D}_f\times\md D_g}\Delta(h;f,g)=\max_{(f,g)\in\mathcal{D}_f\times\md
D_g}\Delta(h^0;f,g).\]
\end{ozn}

Taking into account the introduced definitions and the relations derived in the previous sections  we can verify that the following lemma holds true.

\begin{lema}
The spectral densities $f^0\in\mathcal{D}_f$,
$g^0\in\mathcal{D}_g$ which satisfy the minimality condition (\ref{umova11_f_st.n_d})
are least favorable in the class $\md D=\md D_f\times\md D_g$ for
the optimal linear estimation of the functional $A_N\vec\xi$ based on observations of the sequence $\xi(m)+\eta(m)$
at points   $m\in\mr Z\setminus\{0,1,2,\ldots,N\}$
if the matrices  $(\me T^{\overline{\mu}}_N)^0$, $ (\me P^{\overline{\mu}}_N)^0$, $(\me Q_N)^0$ whose elements are defined by the Fourier coefficients of the functions
\[
\frac{|{\beta^{(d)}(i\lambda)}|^2}{\chi_{\overline{\mu}}^{(d)}(e^{-i\lambda})}
g^0(\lambda)
(f^0(\lambda)+|\beta^{(d)}(i\lambda)|^2 g^0(\lambda))^{-1},\quad
\frac{|{\beta^{(d)}(i\lambda)}|^2}{|\chi_{\overline{\mu}}^{(d)}(e^{-i\lambda})|^2}(f^0(\lambda)+|\beta^{(d)}(i\lambda)|^2 g^0(\lambda))^{-1},
\]
\[
f^0(\lambda)g^0(\lambda)
(f^0(\lambda)+|\beta^{(d)}(i\lambda)|^2 g^0(\lambda))^{-1}
\]
determine a solution of the constrained optimisation problem
\begin{equation*}
   \max_{(f,g)\in \mathcal{D}_f\times\md D_g}\ld(\ld\langle [D_N^{\overline{\mu}}\me a_N]_{+n(\gamma)}- \me T^{\overline{\mu}}_N\me
    a_{\overline{\mu}},(\me P^{\overline{\mu}}_N)^{-1}[D_N^{\overline{\mu}}\me a_N]_{+n(\gamma)}-
   (\me P^{\overline{\mu}}_N)^{-1}\me T^{\overline{\mu}}_N\me a^{\overline{\mu}}_N\rd\rangle\rd.
        +\ld.\langle\me Q_N\me
    a_N,\me a_N\rangle\rd)=
    \end{equation*}
\begin{equation}
    = \ld\langle [D_N^{\overline{\mu}}\me a_N]_{+n(\gamma)}- (\me T^{\overline{\mu}}_N)^0\me
    a^{\overline{\mu}}_N,((\me P^{\overline{\mu}}_N)^0)^{-1}[D_N^{\overline{\mu}}\me a_N]_{+n(\gamma)} -
 ((\me P^{\overline{\mu}}_N)^0)^{-1}(\me T^{\overline{\mu}}_N)^0\me a^{\overline{\mu}}_N
    \rd\rangle
    +\langle\me Q^0_N\me
    a_N,\me a_N\rangle.
\label{minimax1}
\end{equation}
The minimax spectral characteristic $h^0=\vec h_{\overline{\mu},N}(f^0,g^0)$ is calculated by formula (\ref{spectr A}) if
$\vec h_{\overline{\mu},N}(f^0,g^0)\in H_{\mathcal{D}}$.
\end{lema}

The more detailed analysis of properties of the least favorable spectral densities and the minimax-robust spectral characteristics shows that the minimax spectral characteristic $h^0$ and the least favourable spectral densities $f^0$ and $g^0$ form a saddle
point of the function $\Delta(h;f,g)$ on the set
$H_{\mathcal{D}}\times\mathcal{D}$.
The saddle point inequalities
\[
 \Delta(h;f^0,g^0)\geq\Delta(h^0;f^0,g^0)\geq\Delta(h^0;f,g)\quad\forall (f,g)\in
 \mathcal{D},\forall h\in H_{\mathcal{D}}\]
hold true if $h^0=\vec h_{\overline{\mu},N}(f^0,g^0)$,
$\vec h_{\overline{\mu}}(f^0,g^0)\in H_{\mathcal{D}}$ and $(f^0,g^0)$ is a solution of the constrained optimization problem
\be
 \widetilde{\Delta}(f,g)=-\Delta(\vec h_{\overline{\mu}}(f^0,g^0);f,g)\to
 \inf,\quad (f,g)\in \mathcal{D},\label{zad_um_extr_e_d}
 \ee
where the functional $\Delta(\vec h_{\overline{\mu},N}(f^0,g^0);f,g)$ is calculated by the formula

\begin{align*}
\nonumber
&\Delta(\vec h_{\overline{\mu},N}(f^0,g^0);f,g)=
\\
&=\frac{1}{2\pi}\int_{-\pi}^{\pi}
 \frac{|\beta^{(d)}(i\lambda)|^2}{|\chi_{\overline{\mu}}^{(d)}(e^{-i\lambda})|^2}
\left[\chi_{\overline{\mu}}^{(d)}(e^{-i\lambda})(\vec{A}_N(e^{i\lambda}))^{\top}g^0(\lambda) -
(\vec{C}^0_{\mu,N}(e^{i \lambda}))^{\top}
\right]
\times
\\
&\times
(f^0(\lambda)+|\beta^{(d)}(i\lambda)|^2 g^0(\lambda))^{-1}\, f(\lambda)\, (f^0(\lambda)+|\beta^{(d)}(i\lambda)|^2 g^0(\lambda))^{-1}
\times
\\
&\times
\left[\chi_{\overline{\mu}}^{(d)}(e^{-i\lambda}){\vec{A}_N(e^{-i\lambda})}g^0(\lambda) -
\vec{C}^0_{\mu,N}(e^{-i \lambda})
\right]
d\lambda+
\\
&+\frac{1}{2\pi}\int_{-\pi}^{\pi}
\frac{1}{|\chi_{\overline{\mu}}^{(d)}(e^{-i\lambda})|^2}
\left[\chi_{\overline{\mu}}^{(d)}(e^{-i\lambda})(\vec{A}_N(e^{i\lambda}))^{\top}f^0(\lambda) +(\lambda)^{2n}
(\vec{C}^0_{\mu,N}(e^{i \lambda}))^{\top}
\right]
\times
\\
&\times
(f^0(\lambda)+|\beta^{(d)}(i\lambda)|^2 g^0(\lambda))^{-1}\, g(\lambda)\, (f^0(\lambda)+|\beta^{(d)}(i\lambda)|^2 g^0(\lambda))^{-1}
\times
\\
&\times
\left[\chi_{\overline{\mu}}^{(d)}(e^{-i\lambda}) {\vec{A}_N(e^{-i\lambda})} f^0(\lambda) +(\lambda)^{2n}
\vec{C}^0_{\mu,N}(e^{-i \lambda})
\right]
d\lambda,
\end{align*}
where
\[
\vec{C}^0_{\overline{\mu},N}(e^{i \lambda})=\sum_{k=0}^{N+n(\gamma)}
\left(((\me P^{\overline{\mu}}_N)^0)^{-1}[D_N^{\overline{\mu}}\me a_N]_{+n(\gamma)}-((\me P^{\overline{\mu}}_N)^0)^{-1}(\me T^{\overline{\mu}}_N)^0\me a^{\overline{\mu}}_N\right)_k e^{ik\lambda}.
\]

The constrained optimization problem (\ref{zad_um_extr_e_d}) is equivalent to the unconstrained optimization problem
\be \label{zad_unconst_extr_f_st_d}
 \Delta_{\mathcal{D}}(f,g)=\widetilde{\Delta}(f,g)+ \delta(f,g|\mathcal{D})\to\inf,\ee
 where $\delta(f,g|\mathcal{D})$ is the indicator function of the set
$\mathcal{D}$, namely $\delta(f,g|\mathcal{D})=0$ if $(f,g)\in \mathcal{D}$ and $\delta(f,g|\mathcal{D})=+\infty$ if $(f,g)\notin \mathcal{D}$.
The condition
 $0\in\partial\Delta_{\mathcal{D}}(f^0,g^0)$ characterizes a solution $(f^0,g^0)$ of the stated unconstrained optimization problem.
 This condition is the necessary and sufficient condition under which the point $(f^0,g^0)$ belongs to the set of minimums of the convex functional $\Delta_{\mathcal{D}}(f,g)$ \cite{Moklyachuk2015,Rockafellar}.
 Thus, it allows us to find the equations which determine the least favourable spectral densities in some special classes of spectral densities $\md D$.

The form of the functional $\Delta(\vec h_{\overline{\mu}}(f^0,g^0);f,g)$ is suitable for application of the Lagrange method of indefinite
multipliers to the constrained optimization problem \eqref{zad_um_extr_e_d}.
Thus, the complexity of the problem is reduced to  finding the subdifferential of the indicator function of the set of admissible spectral densities. We illustrate the solving of the problem \eqref{zad_unconst_extr_f_st_d} for concrete sets admissible spectral densities  in the following subsections.

\subsection{Least favorable spectral density in classes $\md D_0\times \md D_{\varepsilon}$}
\label{class_D0}

Consider the minimax interpolation problem for the functional $A_N\vec{\xi}$
depending on the unobserved values of the stochastic sequence $\vec{\xi}(m)$ which determine a stationary  GM increments sequence
$\chi_{\overline{\mu},\overline{s}}^{(d)}(\vec{\xi}(m))$  for the following sets of admissible spectral densities $\md D_0^k$, $k=1,2,3,4$,

$$\md D_{0}^{1} =\bigg\{f(\lambda )\left|\frac{1}{2\pi} \int
_{-\pi}^{\pi}
\frac{|\chi_{\overline{\mu}}^{(d)}(e^{-i\lambda})|^2}{|\beta^{(d)}(i\lambda)|^2}
f(\lambda )d\lambda  =P\right.\bigg\},$$
$$\md D_{0}^{2} =\bigg\{f(\lambda )\left|\frac{1}{2\pi }
\int _{-\pi }^{\pi}
\frac{|\chi_{\overline{\mu}}^{(d)}(e^{-i\lambda})|^2}{|\beta^{(d)}(i\lambda)|^2}
{\rm{Tr}}\,[ f(\lambda )]d\lambda =p\right.\bigg\},$$
$$\md D_{0}^{3} =\bigg\{f(\lambda )\left|\frac{1}{2\pi }
\int _{-\pi}^{\pi}
\frac{|\chi_{\overline{\mu}}^{(d)}(e^{-i\lambda})|^2}{|\beta^{(d)}(i\lambda)|^2}
f_{kk} (\lambda )d\lambda =p_{k}, k=\overline{1,T}\right.\bigg\},$$
$$\md D_{0}^{4} =\bigg\{f(\lambda )\left|\frac{1}{2\pi} \int _{-\pi}^{\pi}
\frac{|\chi_{\overline{\mu}}^{(d)}(e^{-i\lambda})|^2}{|\beta^{(d)}(i\lambda)|^2}
\left\langle B_{1} ,f(\lambda )\right\rangle d\lambda  =p\right.\bigg\},$$

\noindent
where  $p, p_k, k=\overline{1,T}$ are given numbers, $P, B_1$ are given positive-definite Hermitian matrices, and sets of admissible spectral densities
$\md D_{\varepsilon}^{k}$, $k=1,2,3,4$ for the stationary noise sequence $\vec \eta(m)$

\begin{equation*}
\md D_{\varepsilon }^{1}  =\bigg\{g(\lambda )\bigg|{\mathrm{Tr}}\,
[g(\lambda )]=(1-\varepsilon ) {\mathrm{Tr}}\,  [g_{1} (\lambda
)]+\varepsilon {\mathrm{Tr}}\,  [W(\lambda )],
\frac{1}{2\pi} \int _{-\pi}^{\pi}
{\mathrm{Tr}}\,
[g(\lambda )]d\lambda =q \bigg\};
\end{equation*}
\begin{equation*}
\md D_{\varepsilon }^{2}  =\bigg\{g(\lambda )\bigg|g_{kk} (\lambda)
=(1-\varepsilon )g_{kk}^{1} (\lambda )+\varepsilon w_{kk}(\lambda),
\frac{1}{2\pi} \int _{-\pi}^{\pi}
g_{kk} (\lambda)d\lambda  =q_{k} , k=\overline{1,T}\bigg\};
\end{equation*}
\begin{equation*}
\md D_{\varepsilon }^{3} =\bigg\{g(\lambda )\bigg|\left\langle B_{2},g(\lambda )\right\rangle =(1-\varepsilon )\left\langle B_{2},g_{1} (\lambda )\right\rangle+\varepsilon \left\langle B_{2},W(\lambda )\right\rangle,
\frac{1}{2\pi}\int _{-\pi}^{\pi}
\left\langle B_{2} ,g(\lambda )\right\rangle d\lambda =q\bigg\};
\end{equation*}
\begin{equation*}
\md D_{\varepsilon }^{4}=\bigg\{g(\lambda )\bigg|g(\lambda)=(1-\varepsilon )g_{1} (\lambda )+\varepsilon W(\lambda ).
\frac{1}{2\pi } \int _{-\pi}^{\pi}
g(\lambda )d\lambda=Q\bigg\},
\end{equation*}

\noindent
where  $g_{1} ( \lambda )$ is a fixed spectral density, $W(\lambda)$ is an unknown spectral density, $q, q_k,k=\overline{1,T}$, are given numbers, $Q$ is a given positive-definite Hermitian matrix.

In the following we will use the next notations:
\begin{multline*}
C^{f0}_{\overline{\mu},N}(e^{i\lambda}):=
\chi_{\overline{\mu}}^{(d)}(e^{i\lambda}){\vec{A}_N(e^{-i\lambda})})^{\top}g^0(\lambda)-
\\
-
\left(\sum_{k=0}^{N+n(\gamma)}
\left(((\me P^{\overline{\mu}}_N)^0)^{-1}[D_N^{\overline{\mu}}\me a_N]_{+n(\gamma)}-((\me P^{\overline{\mu}}_N)^0)^{-1}(\me T^{\overline{\mu}}_N)^0\me a^{\overline{\mu}}_N)_k \right)e^{ik\lambda}\right)^{\top},
\end{multline*}
\begin{multline*}
C^{g0}_{\overline{\mu},N}(e^{i\lambda}):=
\frac{|\chi_{\overline{\mu}}^{(d)}(e^{-i\lambda})|^2}{|\beta^{(d)}(i\lambda)|^2}
({\vec{A}_N(e^{-i\lambda})})^{\top} f^0(\lambda)+
\\
+
\chi_{\overline{\mu}}^{(d)}(e^{-i\lambda})
\left(\sum_{k=0}^{N+n(\gamma)}
\left(((\me P^{\overline{\mu}}_N)^0)^{-1}[D_N^{\overline{\mu}}\me a_N]_{+n(\gamma)}-((\me P^{\overline{\mu}}_N)^0)^{-1}(\me T^{\overline{\mu}}_N)^0\me a^{\overline{\mu}}_N)_k \right)e^{ik\lambda}\right)^{\top},
\end{multline*}
\[
 p_{\chi}^0(\lambda)
=
\frac{|\chi_{\overline{\mu}}^{(d)}(e^{-i\lambda})|^2}{|\beta^{(d)}(i\lambda)|^2} (f^0(\lambda)+{|\beta^{(d)}(i\lambda)|^2}g^0(\lambda)).
\]
From the condition $0\in\partial\Delta_{\mathcal{D}}(f^0,g^0)$
we find the following equations which determine the least favourable spectral densities for these given sets of admissible spectral densities.

For the first pair of the sets of admissible spectral densities $\md D_{f0}^1\times \md D_{\varepsilon} ^{1}$ we have equations
\begin{equation}   \label{eq_4_1}
\left(
C^{f0}_{\overline{\mu},N}(e^{i\lambda})
\right)
\left(
C^{f0}_{\overline{\mu},N}(e^{i\lambda})\right)^{*}
=
 p_{\chi}^0(\lambda)
\vec{\alpha}
\cdot
\vec{\alpha}^{*}
 p_{\chi}^0(\lambda),
\end{equation}

\begin{equation}\label{eq_5_1}
\left(
C^{g0}_{\overline{\mu},N}(e^{i\lambda})
\right)
\left(
C^{g0}_{\overline{\mu},N}(e^{i\lambda})
\right)^{*}
=
(\alpha^{2} +\gamma_1(\lambda ))
\left(
p_{\chi}^0(\lambda)
\right)^2,
\end{equation}

\noindent where $\alpha^{2}$, $\vec{\alpha}$ are Lagrange multipliers,  the function $\gamma_1(\lambda )\le 0$ and $\gamma_1(\lambda )=0$ if ${\mathrm{Tr}}\,[g_{0} (\lambda )]>(1-\varepsilon ) {\mathrm{Tr}}\, [g_{1} (\lambda )]$.

For the second pair of the sets of admissible spectral densities $\md D_{f0}^2\times \md D_{\varepsilon } ^{2}$ we have equation
\begin{equation}  \label{eq_4_2}
\left(
C^{f0}_{\overline{\mu},N}(e^{i\lambda})
\right)
\left(
C^{f0}_{\overline{\mu},N}(e^{i\lambda})\right)^{*}
=
\alpha^{2} \left(
p_{\chi}^0(\lambda)
\right)^{2},
\end{equation}

\begin{equation}   \label{eq_5_2}
\left(
C^{g0}_{\overline{\mu},N}(e^{i\lambda})
\right)
\left(
C^{g0}_{\overline{\mu},N}(e^{i\lambda})
\right)^{*}
=
\left(
p_{\chi}^0(\lambda)
\right)
\left\{(\alpha_{k}^{2} +\gamma_{k}^1 (\lambda ))\delta _{kl} \right\}_{k,l=1}^{T}
\left(
p_{\chi}^0(\lambda)
\right),
\end{equation}

\noindent where $\alpha^{2}$, $\alpha _{k}^{2}$ are Lagrange multipliers,  functions $\gamma_{k}^1(\lambda )\le 0$ and $\gamma_{k}^1 (\lambda )=0$ if $
g_{kk}^{0}(\lambda )>(1-\varepsilon )g_{kk}^{1} (\lambda )$.

For the third pair of the sets of admissible spectral densities $\md D_{f0}^3\times \md D_{\varepsilon }^{3}$ we have equation

\begin{equation}   \label{eq_4_3}
\left(
C^{f0}_{\overline{\mu},N}(e^{i\lambda})
\right)
\left(
C^{f0}_{\overline{\mu},N}(e^{i\lambda})\right)^{*}
=
\left(
p_{\chi}^0(\lambda)
\right)\left\{\alpha _{k}^{2} \delta _{kl} \right\}_{k,l=1}^{T}
\left(
p_{\chi}^0(\lambda)
\right),
\end{equation}

\begin{equation}   \label{eq_5_3}
\left(
C^{g0}_{\overline{\mu},N}(e^{i\lambda})
\right)
\left(
C^{g0}_{\overline{\mu},N}(e^{i\lambda})
\right)^{*}
=
\left(\alpha^{2} +\gamma_1'(\lambda ))
p_{\chi}^0(\lambda)
\right)
B_{2}^{ \top}\,
\left(
p_{\chi}^0(\lambda)
\right),
\end{equation}

\noindent where $\alpha _{k}^{2}$, $\alpha^{2}$ are Lagrange multipliers,  function $\gamma_1' ( \lambda )\le 0$ and $\gamma_1' ( \lambda )=0$ if $\langle B_{2} ,g_{0} ( \lambda ) \rangle>(1- \varepsilon ) \langle B_{2} ,g_{1} ( \lambda ) \rangle$,
 $\delta _{kl}$ are Kronecker symbols.

For the fourth pair of the sets of admissible spectral densities $\md D_{f0}^4\times \md D_{\varepsilon}^{4}$ we have equation
\begin{equation} \label{eq_4_4}
\left(
C^{f0}_{\overline{\mu},N}(e^{i\lambda})
\right)
\left(
C^{f0}_{\overline{\mu},N}(e^{i\lambda})\right)^{*}
=
\alpha^{2} \left(
p_{\chi}^0(\lambda)
\right)
B_{1}^{\top}
 \left(
 p_{\chi}^0(\lambda)
 \right),
\end{equation}

\begin{equation}  \label{eq_5_4}
\left(
C^{g0}_{\overline{\mu},N}(e^{i\lambda})
\right)
\left(
C^{g0}_{\overline{\mu},N}(e^{i\lambda})
\right)^{*}
=
\left(
p_{\chi}^0(\lambda)
\right)
(\vec{\alpha}\cdot \vec{\alpha}^{*}+\Gamma(\lambda))
\left(
p_{\chi}^0(\lambda)
\right),
\end{equation}

\noindent where $\alpha^{2}$, $\vec{\alpha}$ are Lagrange multipliers,  function $\Gamma(\lambda )\le 0$ and $\Gamma(\lambda )=0$ if $g_{0}(\lambda )>(1-\varepsilon )g_{1} (\lambda )$.

The following theorem  holds true.

\begin{thm}
 The least favorable spectral densities $f^0(\lambda)$ and $g^0(\lambda)$ in the classes $ \md  D_0^{k}\times{\md D_{\varepsilon }^{k}}$, $k=1,2,3,4$, for the optimal linear estimation  of the functional  $A_N\vec{\xi}$   are determined by  pairs of equations
\eqref{eq_4_1}-\eqref{eq_5_1}, \eqref{eq_4_2}-\eqref{eq_5_2}, \eqref{eq_4_3}-\eqref{eq_5_3}, \eqref{eq_4_4}-\eqref{eq_5_4},
the minimality condition (\ref{umova11_f_st.n_d}),
the constrained optimization problem (\ref{minimax1}) and restrictions  on densities from the corresponding classes $ \md  D_0^{k}\times{\md D_{\varepsilon }^{k}}$, $k=1,2,3,4$.  The minimax-robust spectral characteristic $\vec h_{\overline{\mu},N}(f^0,g^0)$ of the optimal estimate of the functional $A_N\vec{\xi}$ is determined by the formula (\ref{spectr A}).
\end{thm}

\subsection{Least favorable spectral density in classes $\md D_{1\delta} \times \md D_{V}^{U}$}
\label{class_D0}

Consider the minimax interpolation problem for the functional $A_N\vec{\xi}$
depending on the unobserved values of the stochastic sequence $\vec{\xi}(m)$ which determine a stationary  GM increments sequence
$\chi_{\overline{\mu},\overline{s}}^{(d)}(\vec{\xi}(m))$  for the following sets of admissible spectral densities $\md D_{1\delta}^{k}$, $k=1,2,3,4$,

\begin{equation*}
\md D_{1\delta}^{1}=\left\{f(\lambda )\biggl|\frac{1}{2\pi} \int_{-\pi}^{\pi}
\frac{|\chi_{\overline{\mu}}^{(d)}(e^{-i\lambda})|^2}{|\beta^{(d)}(i\lambda)|^2}
\left|{\rm{Tr}}(f(\lambda )-f_{1} (\lambda))\right|d\lambda \le \delta\right\};
\end{equation*}
\begin{equation*}
\md D_{1\delta}^{2}=\left\{f(\lambda )\biggl|\frac{1}{2\pi } \int_{-\pi}^{\pi}
\frac{|\chi_{\overline{\mu}}^{(d)}(e^{-i\lambda})|^2}{|\beta^{(d)}(i\lambda)|^2}
\left|f_{kk} (\lambda )-f_{kk}^{1} (\lambda)\right|d\lambda  \le \delta_{k}, k=\overline{1,T}\right\};
\end{equation*}
\begin{equation*}
\md D_{1\delta}^{3}=\left\{f(\lambda )\biggl|\frac{1}{2\pi } \int_{-\pi}^{\pi}
\frac{|\chi_{\overline{\mu}}^{(d)}(e^{-i\lambda})|^2}{|\beta^{(d)}(i\lambda)|^2}
\left|\left\langle B_{1} ,f(\lambda )-f_{1}(\lambda )\right\rangle \right|d\lambda  \le \delta\right\};
\end{equation*}
\begin{equation*}
\md D_{1\delta}^{4}=\left\{f(\lambda )\biggl|\frac{1}{2\pi} \int_{-\pi}^{\pi}
\frac{|\chi_{\overline{\mu}}^{(d)}(e^{-i\lambda})|^2}{|\beta^{(d)}(i\lambda)|^2}
\left|f_{ij} (\lambda )-f_{ij}^{1} (\lambda)\right|d\lambda  \le \delta_{i}^j, i,j=\overline{1,T}\right\},
\end{equation*}

\noindent
where  $f_{1} ( \lambda )$ is a fixed spectral density,  $B_1$ is a given positive-definite Hermitian matrix,
$\delta,\delta_{k},k=\overline{1,T}$, $\delta_{i}^{j}, i,j=\overline{1,T}$, are given numbers, and sets of admissible spectral densities $ {\md D_{V}^{U}} ^{k}$, $k=1,2,3,4$ for the stationary noise sequence $\vec \eta(m)$:

\begin{equation*}
 {\md D_{V}^{U}} ^{1}=\left\{g(\lambda )\bigg|V(\lambda )\le g(\lambda
)\le U(\lambda ), \frac{1}{2\pi } \int _{-\pi}^{\pi}
g(\lambda )d\lambda=Q\right\},
\end{equation*}
\begin{equation*}
  {\md D_{V}^{U}} ^{2}  =\bigg\{g(\lambda )\bigg|{\mathrm{Tr}}\, [V(\lambda
)]\le {\mathrm{Tr}}\,[ g(\lambda )]\le {\mathrm{Tr}}\, [U(\lambda )],
\frac{1}{2\pi } \int _{-\pi}^{\pi}
{\mathrm{Tr}}\,  [g(\lambda)]d\lambda  =q \bigg\},
\end{equation*}
\begin{equation*}
{\md D_{V}^{U}} ^{3}  =\bigg\{g(\lambda )\bigg|v_{kk} (\lambda )  \le
g_{kk} (\lambda )\le u_{kk} (\lambda ),
\frac{1}{2\pi} \int _{-\pi}^{\pi}
g_{kk} (\lambda
)d\lambda  =q_{k} , k=\overline{1,T}\bigg\},
\end{equation*}
\begin{equation*}
{\md D_{V}^{U}} ^{4}  =\bigg\{g(\lambda )\bigg|\left\langle B_{2}
,V(\lambda )\right\rangle \le \left\langle B_{2},g(\lambda
)\right\rangle \le \left\langle B_{2} ,U(\lambda)\right\rangle,
\frac{1}{2\pi }
\int _{-\pi}^{\pi}
\left\langle B_{2},g(\lambda)\right\rangle d\lambda  =q\bigg\},
\end{equation*}
where the spectral densities $V( \lambda ),U( \lambda )$ are known and fixed, $ q$,  $q_k$, $k=\overline{1,T}$ are given numbers, $Q$, $B_2$ are given positive definite Hermitian matrices.

From the condition $0\in\partial\Delta_{\mathcal{D}}(f^0,g^0)$
we find the following equations which determine the least favourable spectral densities for these given sets of admissible spectral densities.

For the first pair of the sets of admissible spectral densities $\md D_{1\delta}^{1}  \times {\md D_{V}^{U}} ^{1}$ we have equations

\begin{equation} \label{eq_5_1a}
\left(
C^{f0}_{\overline{\mu},N}(e^{i\lambda})
\right)
\left(
C^{f0}_{\overline{\mu},N}(e^{i\lambda})\right)^{*}
=
\beta^{2} \gamma_2( \lambda )\left(
p_{\chi}^0(\lambda)
\right)^2,
\end{equation}

\begin{equation} \label{eq_5_1b}
 \left(
C^{g0}_{\overline{\mu},N}(e^{i\lambda})
\right)
\left(
C^{g0}_{\overline{\mu},N}(e^{i\lambda})
\right)^{*}
=
\left(
p_{\chi}^0(\lambda)
\right)
(\vec{\beta}\cdot \vec{\beta}^{*}+\Gamma _{1} (\lambda )+\Gamma _{2} (\lambda ))
\left(
p_{\chi}^0(\lambda)
\right),
\end{equation}

\begin{equation} \label{eq_5_1c}
\frac{1}{2 \pi} \int_{-\pi}^{ \pi}
\frac{|\chi_{\overline{\mu}}^{(d)}(e^{-i\lambda})|^2}{|\beta^{(d)}(i\lambda)|^2}
\left|{\mathrm{Tr}}\, (f_0( \lambda )-f_{1}(\lambda )) \right|d\lambda =\delta,
\end{equation}

\noindent
where $\beta^{2}$, $ \vec{\beta}$ are Lagrange multipliers, the function  $\Gamma _{1} (\lambda )\le 0$ and $\Gamma _{1} (\lambda )=0$ if $g^0(\lambda )>V(\lambda )$, the function  $
\Gamma _{2} (\lambda )\ge 0$ and $\Gamma _{2} (\lambda )=0$ if $g^0(\lambda )<U(\lambda )$,
 the function $\left| \gamma_2( \lambda ) \right| \le 1$ and
\[\gamma_2( \lambda )={ \mathrm{sign}}\; ({\mathrm{Tr}}\, (f_{0} ( \lambda )-f_{1} ( \lambda ))): \; {\mathrm{Tr}}\, (f_{0} ( \lambda )-f_{1} ( \lambda )) \ne 0.\]

For the second pair of the sets of admissible spectral densities $\md D_{1\delta}^{2}  \times {\md D_{V}^{U}} ^{2}$ we have equations

\begin{equation}  \label{eq_5_2a}
\left(
C^{f0}_{\overline{\mu},N}(e^{i\lambda})
\right)
\left(
C^{f0}_{\overline{\mu},N}(e^{i\lambda})\right)^{*}
=
\left(
p_{\chi}^0(\lambda)
\right)
\left \{ \beta_{k}^{2} \gamma^2_{k} ( \lambda ) \delta_{kl} \right \}_{k,l=1}^{T}
\left(
p_{\chi}^0(\lambda)
\right),
\end{equation}

\begin{equation}   \label{eq_5_2b}
 \left(
C^{g0}_{\overline{\mu},N}(e^{i\lambda})
\right)
\left(
C^{g0}_{\overline{\mu},N}(e^{i\lambda})
\right)^{*}
=
(\beta^{2} +\gamma _{1} (\lambda )+\gamma _{2} (\lambda ))\left(
p_{\chi}^0(\lambda)
\right)^2,
\end{equation}

\begin{equation} \label{eq_5_2c}
\frac{1}{2 \pi} \int_{- \pi}^{ \pi}
\frac{|\chi_{\overline{\mu}}^{(d)}(e^{-i\lambda})|^2}{|\beta^{(d)}(i\lambda)|^2}
\left|f^0_{kk} ( \lambda)-f_{kk}^{1} ( \lambda ) \right| d\lambda =\delta_{k}, \; k= \overline{1,T},
\end{equation}

\noindent where  $\beta^{2}$, $\beta_{k}^{2}$ are Lagrange multipliers, $\delta _{kl}$ are Kronecker symbols, the function  $\gamma _{1} (\lambda )\le 0$ and $\gamma _{1} (\lambda )=0$ if ${\mathrm{Tr}}\,
[g^0(\lambda )]> {\mathrm{Tr}}\,  [V(\lambda )]$, the function  $\gamma _{2} (\lambda )\ge 0$ and $\gamma _{2} (\lambda )=0$ if $ {\mathrm{Tr}}\,[g^0(\lambda )]< {\mathrm{Tr}}\, [ U(\lambda)]$,
 the functions  $\left| \gamma^2_{k} ( \lambda ) \right| \le 1$ and
\[\gamma_{k}^2( \lambda )={ \mathrm{sign}}\;(f_{kk}^{0}( \lambda)-f_{kk}^{1} ( \lambda )): \; f_{kk}^{0} ( \lambda )-f_{kk}^{1}(\lambda ) \ne 0, \; k= \overline{1,T}.\]

For the third pair of the sets of admissible spectral densities $\md D_{1\delta}^{3}  \times {\md D_{V}^{U}} ^{3}$ we have equations

\begin{equation}   \label{eq_5_3a}
\left(
C^{f0}_{\overline{\mu},N}(e^{i\lambda})
\right)
\left(
C^{f0}_{\overline{\mu},N}(e^{i\lambda})\right)^{*}
=
\beta^{2} \gamma_2'( \lambda )
\left(
p_{\chi}^0(\lambda)
\right)
 B_{1}^{\top}
\left(
p_{\chi}^0(\lambda)
\right),
\end{equation}

\begin{equation}  \label{eq_5_3b}
 \left(
C^{g0}_{\overline{\mu},N}(e^{i\lambda})
\right)
\left(
C^{g0}_{\overline{\mu},N}(e^{i\lambda})
\right)^{*}
=
\left(
p_{\chi}^0(\lambda)
\right)
\left\{(\beta_{k}^{2} +\gamma _{1k} (\lambda )+\gamma _{2k} (\lambda ))\delta _{kl}\right\}_{k,l=1}^{T} \left(
p_{\chi}^0(\lambda)
\right),
\end{equation}

\begin{equation} \label{eq_5_3c}
\frac{1}{2 \pi} \int_{- \pi}^{ \pi}
\frac{|\chi_{\overline{\mu}}^{(d)}(e^{-i\lambda})|^2}{|\beta^{(d)}(i\lambda)|^2}
\left| \left \langle B_{1}, f_0( \lambda )-f_{1} ( \lambda ) \right \rangle \right|d\lambda
= \delta,
\end{equation}

\noindent where $\beta^{2}$, $\beta_{k}^{2}$   are Lagrange multipliers, $\delta _{kl}$ are Kronecker symbols, the function  $\gamma _{1k} (\lambda )\le 0$ and $\gamma _{1k} (\lambda )=0$ if $g_{kk}^{0} (\lambda )>v_{kk} (\lambda )$, the function  $\gamma _{2k} (\lambda )\ge 0$ and $\gamma _{2k} (\lambda )=0$ if $g_{kk}^{0} (\lambda )<u_{kk} (\lambda)$, the function
$\left| \gamma_2' ( \lambda ) \right| \le 1$ and
\[\gamma_2' ( \lambda )={ \mathrm{sign}}\; \left \langle B_{1},f_{0} ( \lambda )-f_{1} ( \lambda ) \right \rangle : \; \left \langle B_{1},f_{0} ( \lambda )-f_{1} ( \lambda ) \right \rangle \ne 0.\]

For the fourth pair of the sets of admissible spectral densities $\md D_{1\delta}^{4}  \times {\md D_{V}^{U}} ^{4}$ we have equations

\begin{equation}  \label{eq_5_4a}
\left(
C^{f0}_{\overline{\mu},N}(e^{i\lambda})
\right)
\left(
C^{f0}_{\overline{\mu},N}(e^{i\lambda})\right)^{*}
=
\left(
p_{\chi}^0(\lambda)
\right)
\left \{ \beta_{ij}( \lambda ) \gamma_{ij} ( \lambda ) \right \}_{i,j=1}^{T}
\left(
p_{\chi}^0(\lambda)
\right),
\end{equation}

\begin{equation} \label{eq_5_4b}
 \left(
C^{g0}_{\overline{\mu},N}(e^{i\lambda})
\right)
\left(
C^{g0}_{\overline{\mu},N}(e^{i\lambda})
\right)^{*}
=
(\beta^{2} +\gamma'_{1}(\lambda )+\gamma'_{2}(\lambda ))\left(
p_{\chi}^0(\lambda)
\right)
B_2^\top  \left(
p_{\chi}^0(\lambda)
\right),
\end{equation}

\begin{equation} \label{eq_5_4c}
\frac{1}{2 \pi} \int_{- \pi}^{ \pi}
\frac{|\chi_{\overline{\mu}}^{(d)}(e^{-i\lambda})|^2}{|\beta^{(d)}(i\lambda)|^2}
\left|f^0_{ij}(\lambda)-f_{ij}^{1}( \lambda ) \right|d\lambda = \delta_{i}^{j},\; i,j= \overline{1,T},
\end{equation}

\noindent where  $\beta^{2}$, $ \beta_{ij}$ are Lagrange multipliers, the function $\gamma'_{1}( \lambda )\le 0$ and $\gamma'_{1} ( \lambda )=0$ if $\langle B_{2},g^0( \lambda) \rangle > \langle B_{2},V( \lambda ) \rangle$, the function  $\gamma'_{2}( \lambda )\ge 0$ and $\gamma'_{2} ( \lambda )=0$ if $\langle
B_{2} ,g^0( \lambda) \rangle < \langle B_{2} ,U( \lambda ) \rangle$,
functions $\left| \gamma_{ij} ( \lambda ) \right| \le 1$ and
\[
\gamma_{ij} ( \lambda )= \frac{f_{ij}^{0} ( \lambda )-f_{ij}^{1} (\lambda )}{ \left|f_{ij}^{0} ( \lambda )-f_{ij}^{1}(\lambda) \right|}: \; f_{ij}^{0} ( \lambda )-f_{ij}^{1} ( \lambda ) \ne 0, \; i,j= \overline{1,T}.
\]

The following theorem  holds true.

\begin{thm}
 The least favorable spectral densities $f^0(\lambda)$ and $g^0(\lambda)$ in the classes $\md D_{1\delta}^{k}  \times {\md D_{V}^{U}} ^{k}$, $k=1,2,3,4$, for the optimal linear estimation  of the functional  $A_N\vec{\xi}$   are determined by  pairs of equations
\eqref{eq_5_1a}-\eqref{eq_5_1c}, \eqref{eq_5_2a}-\eqref{eq_5_2c}, \eqref{eq_5_3a}-\eqref{eq_5_3c}, \eqref{eq_5_4a}-\eqref{eq_5_4c},
the minimality condition (\ref{umova11_f_st.n_d}),
the constrained optimization problem (\ref{minimax1}) and restrictions  on densities from the corresponding classes $\md D_{1\delta}^{k}  \times {\md D_{V}^{U}} ^{k}$, $k=1,2,3,4$,  The minimax-robust spectral characteristic $\vec h_{\overline{\mu},N}(f^0,g^0)$ of the optimal estimate of the functional $A_N\vec{\xi}$ is determined by the formula (\ref{spectr A}).
\end{thm}

\section{Conclusions}

In this article, we present methods of solution of the interpolation problem for stochastic sequences with periodically stationary long memory multiple seasonal increments, or sequences with periodically stationary general multiplicative (GM) increments, introduced in the article by Luz and Moklyachuk \cite{Luz_Mokl_extra_GMI}.
These non-stationary stochastic sequences combine  periodic structure of covariation functions of sequences as well as multiple seasonal factors,
including the integrating one.
A short review of the spectral theory of vector-valued generalized multiple increment sequences is presented.
We describe methods of solution of the interpolation  problem in the case where the spectral densities of the sequence $\xi(m)$ and a noise  sequence $\eta(m)$ are exactly known.
Estimates are obtained by applying the Hilbert space projection technique to the vector sequence $\vec \xi(m)+ \vec \eta(m)$ with stationary GM  increments under the stationary noise sequence $\vec \eta(m)$ uncorrelated with $\vec \xi(m)$.
The case of non-stationary fractional integration is discussed as well.
The minimax-robust approach to interpolation problem is applied in the case of spectral uncertainty where the spectral densities of sequences are not exactly known while, instead, sets of admissible spectral densities are specified.
We propose a representation of the mean square error in the form of a linear functional in $L_1$ space with respect to spectral densities, which allows
us to solve the corresponding constrained optimization problem and describe the minimax (robust) estimates of the functionals.
Relations are described which determine the least favourable spectral densities and the minimax spectral characteristics of the optimal estimates of linear functionals
for a collection of specific classes of admissible spectral densities.
These sets are generalizations of the sets of admissible spectral densities described in a survey article by
Kassam and Poor \cite{Kassam_Poor} for stationary stochastic processes.

\end{document}